# Dynamic analysis of free-free Timoshenko beams on elastic foundation under transverse transient ground deformation

Gersena Banushi

Department of Civil and Environmental Engineering, University of California, Berkeley, Berkeley, CA 94720, USA.

g.banushi@berkeley.edu

**Abstract**

Underground infrastructure, such as pipelines and tunnels, can be vulnerable to the effect of transient ground deformation (TGD) caused by different vibration sources, including earthquakes and traffic loadings. Current design methods are based on simple analytical models that idealize the soil movement as a traveling sinusoidal wave, neglecting both the system inertia and the relative displacement at the soil-structure interface. However, as shown in this paper, this assumption may not be valid for buried large diameter pipelines and tunnels requiring accurate dynamic analysis of the buried structure subjected to TGD in the axial and transverse directions. To efficiently analyse the dynamic response of a buried straight beam subjected to transverse TGD, this study introduces a new semi-analytical model based on the Timoshenko beam on Winkler foundation theory. The closed-form analytical solution of the governing differential equation revealed that the vibration spectrum is divided in four parts, separated by three transition frequencies. Across each transition frequency, the oscillatory characteristics of the vibration modes change as a function of the system's mechanical and geometric properties, considerably affecting the dynamic amplification response of the buried Timoshenko beam under TGD. To verify the validity of the proposed model, this work analyses the case study of a buried 107 cm (42 inch) diameter steel water pipeline of varying lengths and operating conditions, subjected to transverse TGD. Comparison of the obtained analytical solutions with the finite element analysis results showed excellent agreement between the two approaches, demonstrating the accuracy of the proposed model. Additional validation was achieved through comparison with high-resolution experimental data from free-vibration tests and multi-point shaking-table experiments, further confirming the accuracy and robustness of the proposed model. The frequency response analysis revealed dynamic amplification of the soil-structure interaction for forcing frequencies near the system's fundamental frequency. These may fall within the range of dominant frequencies characterizing seismic waves, requiring accurate dynamic analysis. The proposed methodology provides a robust analytical framework for evaluating the primary factors impacting the dynamic behavior of buried beams, giving a deeper understanding of the system response under various sources of ground vibration.

## Introduction

Earthquakes and traffic loadings generate waves in the ground propagating at various frequencies and apparent velocities and the resulting transient ground deformations (TGDs) potentially affect a large part of underground infrastructure systems such as pipelines and tunnels. For example, seismic wave propagation has caused significant damage to buried lifelines in the 1985 Michoacan earthquake event. A damage ratio of about 0.45 repairs/km has been reported for the water supply system in the soft soil zone of Metropolitan Mexico City (O'Rourke and Liu, 2012). While high frequency seismic ground motions present dominant frequencies on the range 5-20 Hz (McNutt, 1996; Abbasiverki and Ansell, 2020), traffic induced ground vibration falls in the range of 10-30 Hz (Hendriks, 2002). This paper discusses that accurate analysis of the dynamic behaviour of buried infrastructure may be required for certain scenarios under these frequency ranges, and this is carried out using a newly introduced semi-analytical model.

Existing dynamic analysis methods of underground lifelines are broadly classified in two broad groups (Manolis and Beskos, 1997): (i) experimental and real site observations, and (ii) analytical and numerical procedures.
Experimental studies, including free-vibration testing of beams, shaking-table simulations, centrifuge experiments, and full-scale site measurements, have significantly advanced understanding of beam vibration and dynamic soil-structure interaction mechanisms (Monsivais et al., 2016; Han et al., 2024). However, experimental investigations are often limited to specific geometries, boundary conditions, material types, and thus serve primarily as benchmarks for validating analytical and numerical models.

A key design factor is the soil-structure interaction, that depends on the stiffness difference between the structure and the surrounding soil. If this is neglected, the structure follows the free field ground motion, and the strains at either side of the soil-structure interface are nearly the same. Current seismic design analysis for buried pipelines (ALA 2001; ASCE 1984; PRCI 2004; CEN, 2003; IITK-GSDMA. 2007) is based on Newmark's analytical model (Newmark, 1976) that idealizes the soil movement as a traveling sinusoidal wave, while ignoring the system's inertia and the relative movement at the soil-structure interface. This results in pipeline displacements and strains that are equal to those of the ground. Kuesel (1969) and Yeh (1974) extended the Newmark's solution to account for S- and Rayleigh waves propagating at an oblique angle relative to the structure axis. The axial and bending strains along the pipeline are evaluated by assessing separately each direction of the seismic induced TGD, i.e., parallel and perpendicular to the pipeline axis. The resulting pipeline bending strain associated with the pipeline curvature are considered a second order effect compared to the axial strains, and typically neglected (ASCE, 1984; O'Rourke and Liu, 2012). However, this assumption may not be valid for buried large diameter pipelines and tunnels and accurate dynamic analysis of the buried structure subjected to TGD in the two horizontal directions is required (Zerva 1993, Ariman and Muleski, 1981; Anastasopoulos et al., 2007). This study shows that dynamic application effects may need to be considered for cases with poorly compacted soils with small soil stiffness and/or the increased effective mass of fluid filled pipeline causing a reduction in the natural frequencies.

Advanced analytical models are available for this problem by assuming the structure as a cylindrical shell embedded in an elastic half space (Datta et al., 1984; Liu et al., 1991) or by using the beam on Winkler foundation model. The latter approach has been extensively adopted in different lifeline engineering applications, including railway tracks, tunnels and buried pipelines. Within this method, the system has been analyzed either quasi-statically, i.e. by ignoring the inertia effects (Takada and Tanabe, 1987; O'Rourke and Wang, 1978; Kiyomiya, 1995; O'Rourke et al., 2004; Wang et al., 2005; McLaughlin and O'Rourke, 2009; Banushi and Weidlich, 2018; O'Rourke and Liu, 2012), or dynamically (Sakurai and Takahashi, 1969; Xu et al., 2021).

Various numerical procedures, including the finite element method (FEM), have been extensively employed to evaluate the dynamic response of buried beams, allowing to analyze more complex system configurations in structural design practice and research (St John and Zahrah, 1987; Yang et al., 1988; Jahangiri and Shakib, 2018). Numerical models have often been used to verify the validity of proposed analytical and semi-analytical solutions of this soil-structure interaction problem (Kouretzis et al., 2006; Morfidis, 2010; Wu, 2019). These models also represent a useful tool to validate more advanced numerical models (Anastasopoulos et al., 2007) and can be more easily integrated in the design guidelines of buried infrastructure.

The dynamic analysis of buried beams includes the free structural vibration and the forced transient response, using either the Euler-Bernoulli or the Timoshenko beam theory. The former approach has been more widely employed in dynamic analysis of a structure resting on Winkler foundation, including railways, tunnels and pipelines, because of the simpler mathematical formulation (Wang, 1978; Mavridis and Pitilakis, 1996; Hindy and Novak, 1980; Zerva, 1993; Hosseini and Roudsari, 2015). These solutions consider different beam lengths and boundary conditions. The Bernoulli-Euler theory is not valid for the cases of non slender beams and free-free or pinned-free shear beams (Kausel, 2002; Balkaya et al., 2009). On the other hand, the Timoshenko model considers the effect of transverse shear deformation and rotary inertia on the beam modes and frequencies, providing a more accurate dynamic response (Majkut, 2009). Among various boundary conditions, the free-free Timoshenko beam on elastic foundation has been less investigated, either analytically or numerically, because of the more complex computation, including the rigid-body motion (Kuessel, 2002; Wu, 2019).
Free–free boundary conditions are appropriate for pipeline or tunnel segments connected through flexible joints such as expansion couplings or ball-and-socket joints (ALA, 2005; Sun et al., 2017). They also represent unconstrained sections, including those near a right-angled bend at the base of a riser or at pipeline terminations where the pipe enters a tank and the annular penetration gap is filled with soft material, providing negligible restraint (O'Rourke et al., 2008).
Moreover, in longer pipelines, the undamped natural frequencies of all lateral modes become closer to each other, requiring a significantly greater number of modes to calculate dynamic response (Hindy and Novak, 1980; Zerva et al., 1988; Hosseini and Roudsari, 2015) and consequent computational effort.

To analyse the dynamic response of a buried straight beams with free ends subjected to transverse TGD, this study develops a new semi-analytical model, based on the Timoshenko beam on Winkler foundation theory. The Fourier method of variable separation is used to solve the equation of vibration of the buried Timoshenko beam. The model can evaluate the free vibration and forced dynamic response for various system lengths and operating conditions. The application of the model is illustrated by simulating the case of a buried straight water pipeline with variable lengths and effective masses.
Its accuracy is verified through finite element modal analysis and implicit dynamic simulations, and further validated using high-resolution experimental data on free-vibration behaviour (Monsivais et al., 2016) and multi-point shaking-table tests simulating transverse ground excitation (Han et al., 2024). Finally, the model is used to examine the influence of key system parameters on dynamic response and to characterize dynamic amplification mechanisms.

## 2. Methodology

This section describes the methodology for evaluating the dynamic response of a Timoshenko beam on Winkler foundation through a semi-analytical modal formulation. First, the governing differential equation is solved in a closed form using the Fourier method of variable separation. Second, the transition frequencies $\tilde{\omega}_i$, which partition the vibration spectrum into distinct oscillatory regimes, are evaluated as functions of the relative stiffness of the beam–foundation system and the effects of shear deformation and rotary inertia. Third, the natural vibration frequencies and corresponding mode shapes are obtained over the entire spectrum, considering free–free boundary conditions. The complete free and forced vibration is then calculated explicitly using the representative case study of a Timoshenko beam on Winkler foundation modelling a buried pipeline subjected to seismic wave propagation, as schematically illustrated in Figure 1.

The accuracy of the derived semi-analytical solution is assessed through two complementary validation strategies combining finite-element analysis and laboratory measurements. Specifically, the analytically predicted frequencies and modes for various beam lengths $L$ are compared with the results of a modal analysis performed using ABAQUS/Standard (Simulia 2022). The forced dynamic response is evaluated by numerically solving the modal equations in the unknown modal coordinates and reconstructing the total displacement as a linear combination of all modal contributions. The resulting responses are compared with the finite element analysis predictions, demonstrating the validity of the developed semi-analytical model.

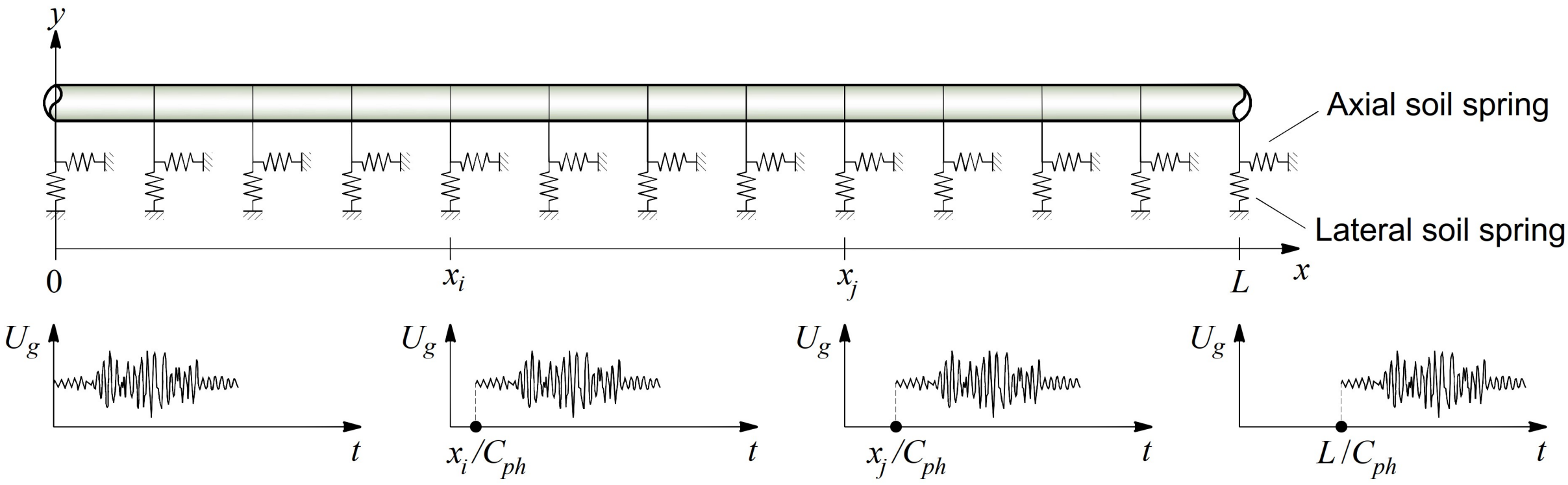

Figure 1. Continuous buried pipeline model subjected to seismic induced ground displacement time histories $U_g(x,t)$, with a time shift at axial distance $x_i$ equal to $x_i/C_{ph}$, inversely proportional to the apparent wave velocity $C_{ph}$.

Furthermore, the proposed formulation is validated using experimental measurements of the free-vibration behaviour of Timoshenko beams (Monsivais et al., 2016) and multi-point shaking-table tests simulating transverse ground excitation on buried pipelines (Han et al., 2024). Finally, the implemented model is used to perform a frequency-response analysis for various beam lengths and operating conditions, providing deeper insight into the system's dynamic behaviour. A flowchart summarizing the complete semi-analytical framework is presented in Figure 2.

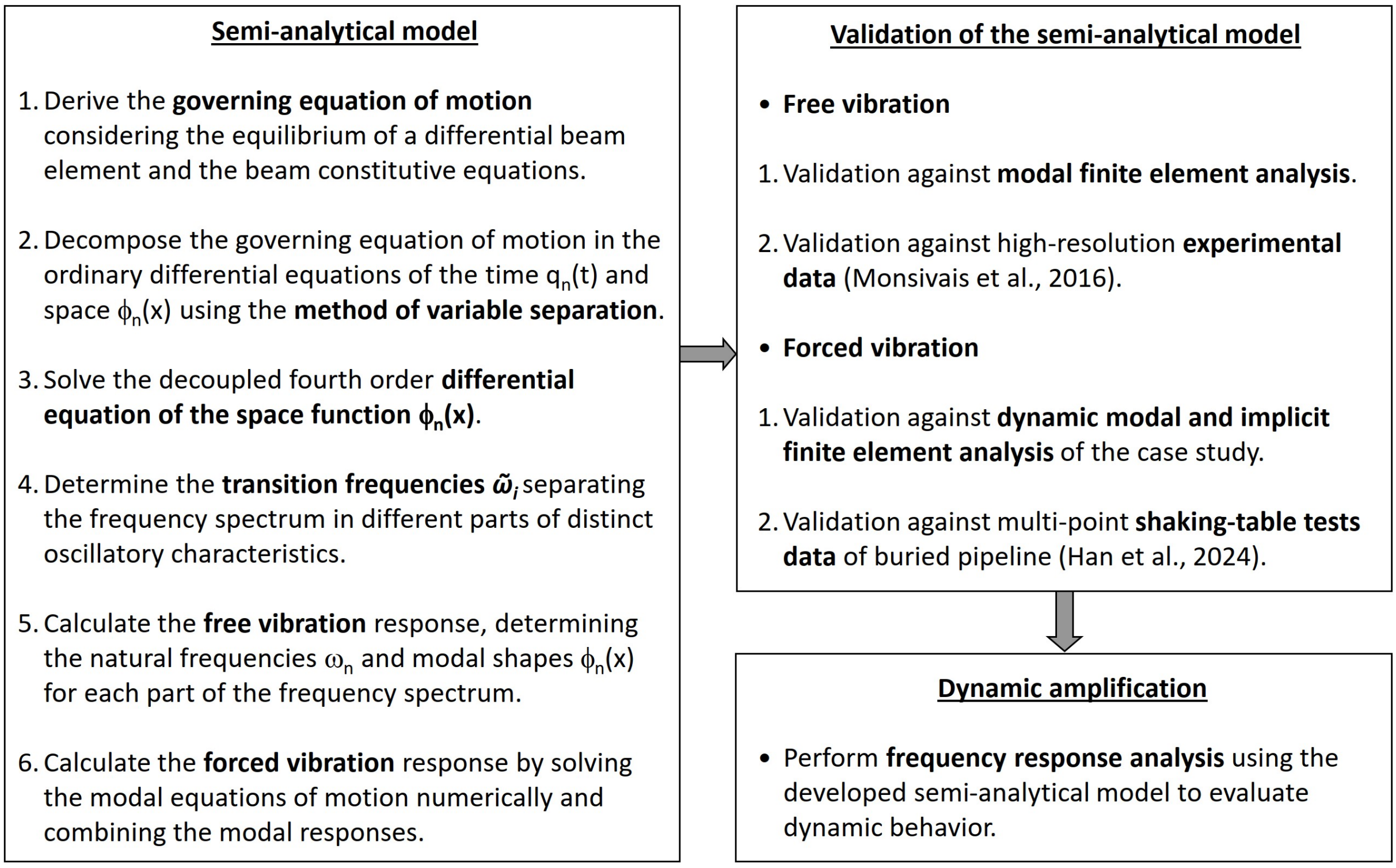

Figure 2. Flowchart of the proposed methodology for assessing the dynamic behavior of buried Timoshenko beam.

## 3. Semi-analytical model of the buried Timoshenko beam under transient ground deformation (TGD)

This section develops the semi-analytical framework for evaluating the dynamic response of a Timoshenko beam on Winkler foundation subjected to transient ground deformation (TGD). First, the governing differential equation of motion is solved in closed form using the Fourier method of variable separation. The transition frequencies $\tilde{\omega}_i$, which partition the vibration spectrum into distinct oscillatory regimes, are then obtained as functions of the relative stiffness of the beam–foundation system, shear deformation, and rotary inertia. Subsequently, the natural vibration frequencies and associated mode shapes are derived for the entire spectrum under free–free boundary conditions. Finally, the modal equations governing the forced vibration response are formulated, enabling numerical computation of the total displacement response.

### 3.1. Exact solution of the equation of vibration of a Timoshenko beam on Winkler foundation

The equation of vibration of a Timoshenko beam in the transverse direction is obtained considering translational and rotational equilibrium of a differential beam element subjected to lateral excitation resting on elastic foundation (Eq. (1)), as illustrated in Figure 3.

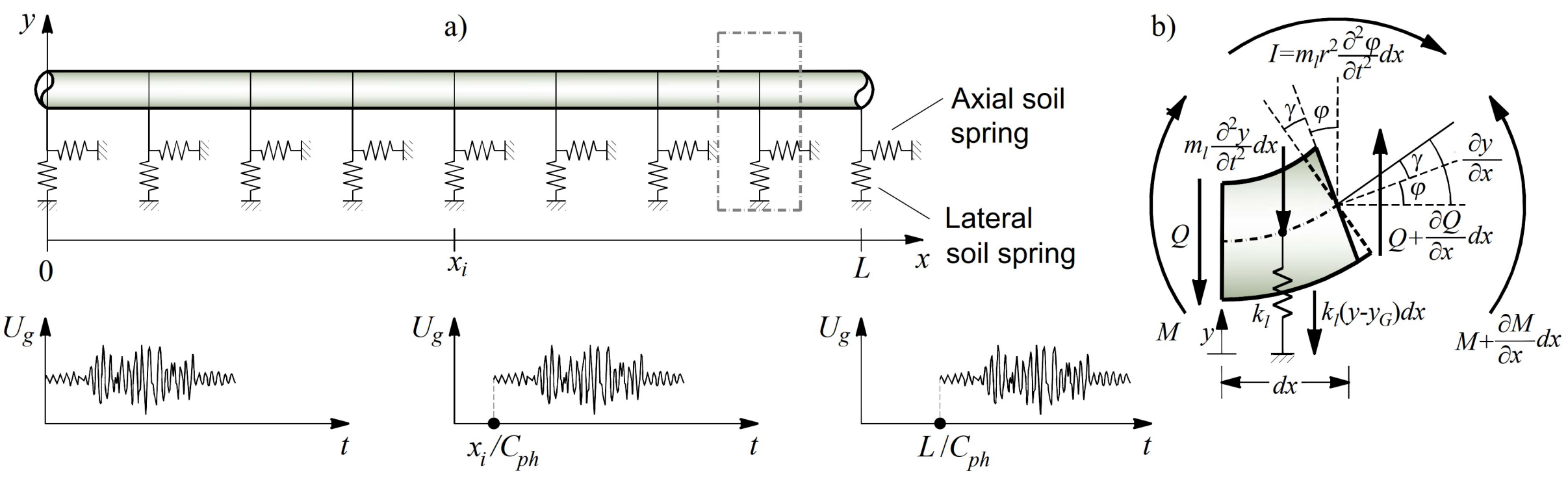

Figure 3 Schematic representation of the a) Timoshenko beam on Winkler foundation subjected to ground displacement $U_g(x,t)$; b) translational and rotational equilibrium of a differential beam element of length $dx$ under lateral deformation, where the soil's resistance to the beam lateral movement is modeled by a linear spring with stiffness $k_l$.

$$\begin{cases} -\dfrac{\partial Q}{\partial x} + m_l \dfrac{\partial^2 y}{\partial t^2} + k_l y = k_l y_G \\ \dfrac{\partial M}{\partial x} + Q - m_l r^2 \dfrac{\partial^2 \varphi}{\partial t^2} = 0 \end{cases} \tag{1}$$

where $Q$ and $M$ are the beam shear force and bending moment, respectively, $m_l$ is the linear mass of the beam per unit length, $m_l$ is the beam linear mass, $k_l$ is the stiffness of the lateral soil spring, $r$ is the radius of gyration of the beam cross section, $y_G$ is the ground displacement in the transverse direction, while $y$ and $\varphi$ are the centroid transverse displacement and cross-sectional rotation, respectively.

By substituting the beam constitutive equations in the system of equations, Eq. (1), we obtain the coupled system of the equations of motion, in terms of the kinematic variables $y$ and $\varphi$:

$$\begin{cases} -kGA_b\left(\dfrac{\partial^2 y}{\partial x^2} - \dfrac{\partial \varphi}{\partial x}\right) + m_l \dfrac{\partial^2 y}{\partial t^2} + k_l y = k_l y_G \\ EJ\dfrac{\partial^2 \varphi}{\partial x^2} + kGA_b\left(\dfrac{\partial y}{\partial x} - \varphi\right) - m_l r^2 \dfrac{\partial^2 \varphi}{\partial t^2} = 0 \end{cases} \tag{2}$$

where $EJ$ is the beam flexural rigidity, $G$ is the modulus of rigidity, $A_b$ is the beam cross sectional area, and $k$ is a constant that depends on the cross-sectional shape and accounts for the nonuniform distribution of shear stress across the cross-section.

This system of the two second-order partial differential equations (PDEs) in Eq. (2) can be reduced to a unique fourth-order PDE, after some algebraic manipulations:

$$EJ\frac{\partial^4 y}{\partial x^4} - m_l r^2\left(1+\frac{E}{kG}\right)\frac{\partial^4 y}{\partial x^2 \partial t^2} + m_l\frac{\partial^2 y}{\partial t^2} + \frac{m_l^2 r^2}{kGA_b}\frac{\partial^4 y}{\partial t^4} + k_l\left(y - \frac{EJ}{kGA_b}\frac{\partial^2 y}{\partial x^2} + \frac{m_l r^2}{kGA_b}\frac{\partial^2 y}{\partial t^2}\right)$$
$$= k_l\left(y_G - \frac{EJ}{kGA_b}\frac{\partial^2 y_G}{\partial x^2} + \frac{m_l r^2}{kGA_b}\frac{\partial^2 y_G}{\partial t^2}\right) \quad (3)$$

The Fourier method of variable separation is used to express the beam displacement and rotation as a linear combination of the modal displacements $y_n(x,t)$ and rotations $\varphi_n(x,t)$, respectively:

$$y(x,t) = \sum_{n=1}^{\infty} y_n(x,t) = \sum_{n=1}^{\infty} \phi_n(x) q_n(t), \quad \varphi(x,t) = \sum_{n=1}^{\infty} \varphi_n(x,t) = \sum_{n=1}^{\infty} \varphi_n(x) q_n(t) \quad (4)$$

where, $q_n(t)$, $\phi_n(x)$, and $\varphi_n(x)$ are the modal coordinate, modal shape and rotation corresponding to the $n$-th natural vibration frequency $\omega_n$, respectively.

Substituting Eq. (4) into the expression of the free vibration of the Timoshenko beam on elastic foundation (left-hand side of Eq. (3)), separates the equation into two ordinary differential equations: one for the time function $q_n(t)$ and one for the spatial function $\phi_n(x)$:

$$\ddot{q}_n(t) + \omega_n^2 q_n(t) = 0 \quad (5.1)$$
$$\frac{EJ}{m_l}\phi_n^{IV}(x) + r^2\left[\omega_n^2\left(1+\frac{E}{kG}\right) - \frac{k_l}{m_l}\frac{E}{kG}\right]\phi_n^{II}(x) + \left[\frac{m_l r^2}{kGA_b}\omega_n^4 - \left(1+\frac{k_l r^2}{kGA_b}\right)\omega_n^2 + \frac{k_l}{m_l}\right]\phi_n(x) = 0 \quad (5.2)$$

or,

$$A\phi_n^{IV}(x) + B\phi_n^{II}(x) + C = 0 \quad (6)$$

where:

$$A = \frac{EJ}{m_l} \quad (6.1)$$
$$B = r^2\left[\omega_n^2\left(1+\frac{E}{kG}\right) - \frac{k_l}{m_l}\frac{E}{kG}\right] \quad (6.2)$$
$$C = \frac{m_l r^2}{kGA_b}\omega_n^4 - \left(1+\frac{k_l r^2}{kGA_b}\right)\omega_n^2 + \frac{k_l}{m_l} = \frac{m_l r^2}{kGA_b}\left(\omega_n^2 - \frac{k_l}{m_l}\right)\left(\omega_n^2 - \frac{kGA_b}{m_l r^2}\right) \quad (6.3)$$

The solutions of the fourth-order ordinary differential equation Eq. (6) have the form of exponential functions $\exp(\lambda x)$, where, in general, $\lambda \in \boldsymbol{C}$. The characteristic equation associated with Eq. (6) is the biquadratic equation:

$$A\lambda^4 + B\lambda^2 + C = 0 \quad (7)$$

having the following solutions:

$$\lambda_1^2 = \frac{-B+\sqrt{\Delta}}{2A} \tag{8}$$

$$\lambda_2^2 = \frac{-B-\sqrt{\Delta}}{2A} \tag{9}$$

where Δ represents the discriminant of the biquadratic Eq. (6):

$$\Delta = B^2 - 4AC = a_\Delta \omega_n^4 + b_\Delta \omega_n^2 + c_\Delta \tag{10}$$

with,

$$a_\Delta = r^4\left(1-\frac{E}{kG}\right)^2 > 0 \tag{11.1}$$

$$b_\Delta = \frac{k_l}{m_l}\left[2\left(1-\frac{E}{kG}\right)\frac{E}{kG}r^4 + \frac{4EJ}{k_l}\right] \tag{11.2}$$

$$c_\Delta = \left(\frac{k_l}{m_l}\right)^2\left[\left(\frac{E}{kG}\right)^2 r^4 - \frac{4EJ}{k_l}\right] \tag{11.3}$$

The term $4EJ/k_l$ in Eqs. (11.2) and (11.3) is equal to the fourth power of the characteristic length of the system, $l_c$ (Hetényi and Hetbenyi, 1946), representing the flexibility of the buried beam relative to the soil stiffness. This value is relatively large for beams resting on soft foundations, determining the sign of the coefficients, $b_\Delta > 0$ and $c_\Delta < 0$ in the Eqs. (11.2) and (11.3), respectively.
From Eq. (10), it results that the discriminant $\Delta$ is positive for values of vibration frequency $\omega_n$ greater than $\tilde{\omega}_1$:

$$\tilde{\omega}_1^2 = \frac{-b_\Delta + \sqrt{b_\Delta^2 - 4a_\Delta c_\Delta}}{2a_\Delta} = \frac{k_l}{m_l}\frac{2\left(\frac{E}{kG}-1\right)\frac{E}{kG} - \left(\frac{l_c}{r}\right)^4 + \sqrt{\left(\frac{l_c}{r}\right)^8 - 4\left(\frac{E}{kG}-1\right)\left(\frac{l_c}{r}\right)^4}}{2\left(\frac{E}{kG}-1\right)^2} \tag{12}$$

Moreover, the roots $\pm\sqrt{\lambda_1^2}$ given by Eq. (8) are real numbers for negative values of C (Eq. 6.3), corresponding to values of vibration frequency $\omega_n$ in the range between $\tilde{\omega}_2$ and $\tilde{\omega}_3$:

$$\tilde{\omega}_2^2 = \frac{k_l}{m_l} \tag{13}$$

$$\tilde{\omega}_3^2 = \frac{kGA_b}{m_l r^2} = \frac{1}{4}\frac{k_l}{m_l}\left(\frac{E}{kG}\right)^{-1}\left(\frac{l_c}{r}\right)^4 \tag{14}$$

whereas the roots $\pm\sqrt{\lambda_2^2}$ given by Eq. (9) are purely imaginary conjugate numbers, for positive values of the discriminant $\Delta$ ($\omega > \tilde{\omega}_1$).
The transition frequency $\tilde{\omega}_3$ (Eq. (14)) represents the classical cut-off frequency for a bare Timoshenko beam ($k_l$ = 0), characterized by the effects of shear deformation and rotary inertia, consistent with observations in previous studies (Majkut, 2009; Monsivais et al., 2016). When the beam rests on an elastic foundation ($k_l \neq 0$), an additional cut-off frequency $\tilde{\omega}_2$ (Eq. (13)) arises due to the contribution of the foundation stiffness. The lowest transition frequency $\tilde{\omega}_1$ (Eq. (12)) depends on the combined influence of the beam shear deformation, rotary inertia, and the flexibility of the buried beam relative to the surrounding soil. In the special case of an Euler Bernoulli beam on a Winkler foundation, where shear deformation and rotary inertia are absent, only the single transition frequency given by Eq. (13) remains.
According to Eq. (12), the ratio between the first two transition frequencies $\tilde{\omega}_1/\tilde{\omega}_2$ is controlled by the ratio between the characteristic system length $l_c = \sqrt[4]{4EJ/k_l}$ and the radius of gyration $r$, as well as the factor $E/kG$. As $l_c/r$ increases, the frequency ratio ($\tilde{\omega}_1/\tilde{\omega}_2$) rises monotonically from its minimum value $(\tilde{\omega}_1/\tilde{\omega}_2)_{\min}$ (Eq. (15)) and asymptotically approaches unity, as shown in Figure 4.

$$\left(\frac{\tilde{\omega}_1}{\tilde{\omega}_2}\right)_{\min} = \sqrt{1-\frac{1}{\dfrac{E}{kG}-1}}, \quad \text{for } \frac{l_c}{r} = \sqrt[4]{4\left(\frac{E}{kG}-1\right)} \tag{15}$$

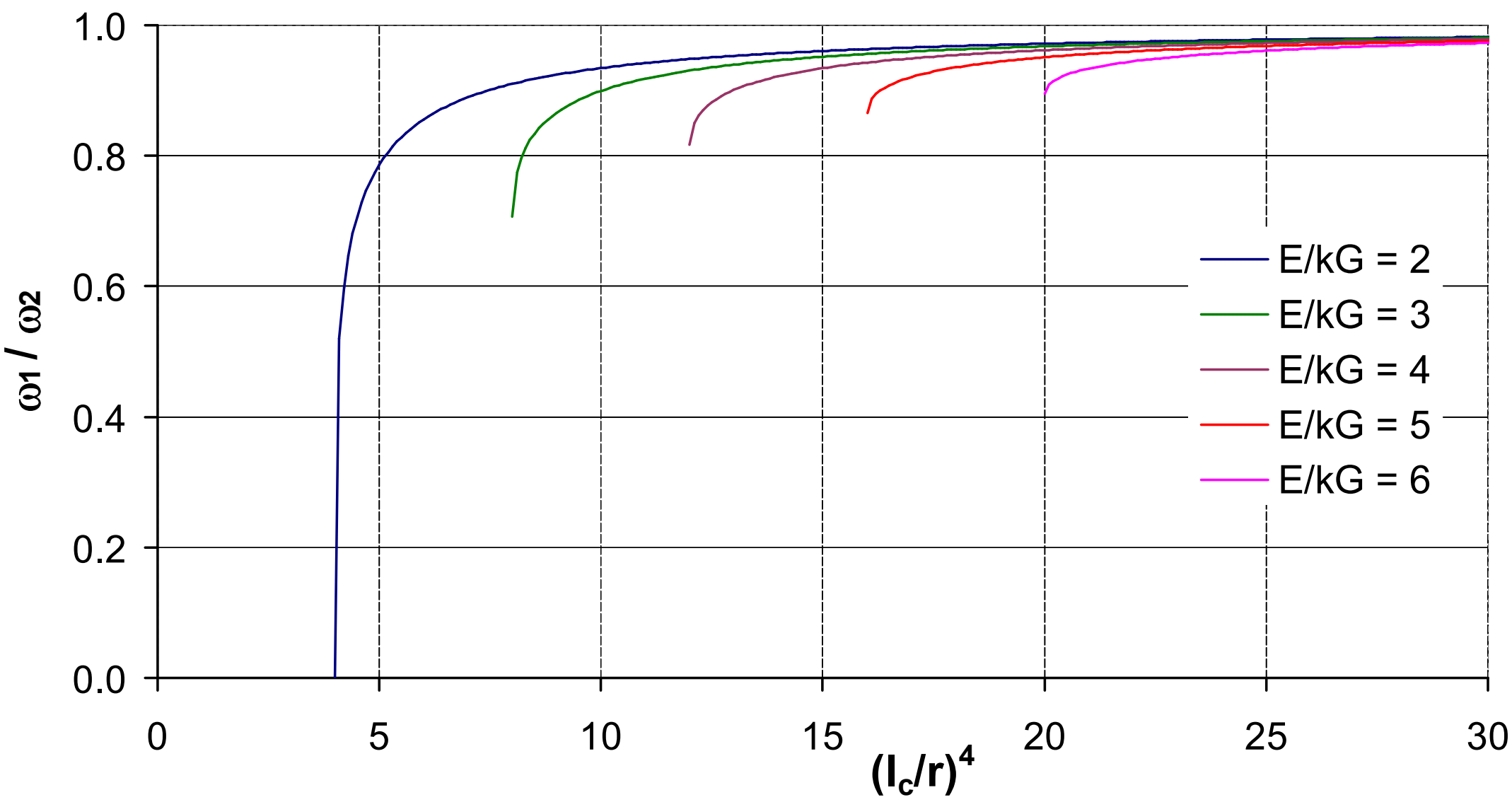

Figure 4. Variation of the frequency ratio $\tilde{\omega}_1/\tilde{\omega}_2$ as a function the dimensionless ratio $l_c/r$ for different values of $E/kG$.

Clearly, the minimum frequency ratio decreases for smaller values of $E/kG$, approaching zero as $E/kG \to 2$. For a Timoshenko beam with tubular cross section ($k = 0.53$) and elastic material with Possion's ratio $\nu = 0.3$, the minimum frequency ratio is $(\tilde{\omega}_1/\tilde{\omega}_2)_{\min} = 0.88$, occurring at $l_c/r = 3.2$. For beams resting on soft foundations, the relative soil-structure flexibility becomes large ($l_c/r >> 1$), and $\tilde{\omega}_1$ approaches $\tilde{\omega}_2$, in agreement with the findings of Banushi and Soga (2025). Conversely, the frequency ratio $\tilde{\omega}_3/\tilde{\omega}_2$ increases with the fourth power of the non-dimensional ratio $l_c/r$, and is inversely proportional to the factor $E/kG$ (Eq. (14)).
The three transition frequencies $\tilde{\omega}_1$, $\tilde{\omega}_2$, $\tilde{\omega}_3$ mark the boundaries between distinct solution regimes of the governing differential Eq. (6). Whether these transition values are part of the frequency spectrum depends on the applied boundary conditions. Across the transition frequency, the oscillatory characteristics of the vibration modes change, significantly affecting the dynamic response of the Timoshenko beam on Winkler foundation, as demonstrated in this study.

### 3.2. Natural vibration frequencies and modes of the Timoshenko beam on Winkler foundation

Considering the assumed values of the system parameters in the present study, $\tilde{\omega}_1 < \tilde{\omega}_2 < \tilde{\omega}_3$. Therefore, the eigenfrequencies $\omega_n$, and associated modal shapes $\phi_n(x)$, can be divided into five cases:

1. $\omega_n < \tilde{\omega}_1$:

The roots of the characteristic Eq. (7) are complex conjugates:

$$\lambda = \pm\sqrt{\frac{-B \pm \sqrt{|\Delta|}\cdot i}{2A}} = \pm\alpha \pm \beta i \tag{16}$$

Where:

$$\alpha = \sqrt{\frac{\sqrt{\left(\frac{-B}{2A}\right)^2 + \left(\frac{\sqrt{|\Delta|}}{2A}\right)^2} + \left(\frac{-B}{2A}\right)}{2}} \tag{17.1}$$

$$\beta = \sqrt{\frac{\sqrt{\left(\frac{-B}{2A}\right)^2 + \left(\frac{\sqrt{|\Delta|}}{2A}\right)^2} - \left(\frac{-B}{2A}\right)}{2}} \tag{17.2}$$

leading to the solution:

$$\phi_n(x) = e^{-\alpha x}\left(C_1 \sin\beta x + C_2 \cos\beta x\right) + e^{\alpha x}\left(C_3 \sin\beta x + C_4 \cos\beta x\right) \tag{18}$$

2. $\omega_n = \tilde{\omega}_1$:

The two roots $\lambda_1^2$ and $\lambda_2^2$ of the characteristic Eq. (7) do coincide. Since $B > 0$ (i.e. $\tilde{\omega}_1^2 > 1/(1+kG/E)\cdot k_l/m_l$), the solution consists of a couple of purely imaginary conjugate roots, $\pm i\lambda$, whose multiplicity is two:

$$\lambda_1^2 = \lambda_2^2 = \left.\frac{-B}{2A}\right|_{\omega=\tilde{\omega}_1} = i^2\lambda^2 \tag{19}$$

leading to the mode shape solution:

$$\phi_n(x) = C_1 \sin\lambda x + C_2 \cos\lambda x + x\left(C_3 \sin\lambda x + C_4 \cos\lambda x\right) \tag{20}$$

3. $\tilde{\omega}_1 < \omega_n < \tilde{\omega}_2$ and $\omega_n > \tilde{\omega}_3$:

Considering Eqs (8) and (9), there are a couple of purely imaginary conjugate roots, i. e. $\pm i\sqrt{\lambda_1^2}$ and $\pm i\sqrt{\lambda_2^2}$, leading to the mode shape solution:

$$\phi_n(x) = C_1 \sin\lambda_1 x + C_2 \cos\lambda_1 x + C_3 \sin\lambda_2 x + C_4 \cos\lambda_2 x \tag{21}$$

4. $\omega_n = \tilde{\omega}_2$ and $\omega_n = \tilde{\omega}_3$:

In this case, there is a zero real root whose multiplicity is two $\pm\sqrt{\lambda_1^2} = 0$, and one couple of purely imaginary roots, $\pm i\sqrt{\lambda_2^2} = \pm i\sqrt{B/A}$, leading to the mode shape solution:

$$\phi_n(x) = C_1 + C_2 x + C_3 \sin\lambda_2 x + C_4 \cos\lambda_2 x \tag{22}$$

5. $\tilde{\omega}_2 < \omega_n < \tilde{\omega}_3$:

For this range of vibration frequency, there are two real roots, $\pm\sqrt{\lambda_1^2}$ and two purely imaginary conjugate roots $\pm i\sqrt{\lambda_2^2}$, leading to the modal shape solution:

$$\phi_n(x) = C_1 \sinh\lambda_1 x + C_2 \cosh\lambda_1 x + C_3 \sin\lambda_2 x + C_4 \cos\lambda_2 x \tag{23}$$

These solutions contain four unknown constants $C_1$, $C_2$, $C_3$, $C_4$, and the eigenvalue parameters $\lambda_1$ and $\lambda_2$. By applying the four end boundary conditions for a single-span beam provides a solution for the natural frequency $\omega_n$, and for the three constants in terms of the fourth, resulting in the natural mode shapes $\phi_n(x)$.

### 3.3. Boundary conditions

The whole frequency spectrum and associated modal shapes depend on the applied boundary conditions. This study evaluates the spectrum for the single span Timoshenko beam of length $L$, considering free end constraints.

The corresponding constrained variables including the shear $Q$, bending moment $M$ and transverse displacement $y$, as well as the equivalent differential equations in terms of the modal shapes $\phi_n(x)$, are indicated in Table 1, where the coefficients $K_M$, and $K_T$ are given by:

$$K_M = \frac{m_l \omega_n^2 - k_l}{kGA_b} \tag{24}$$

$$K_T = K_M + \frac{m_l r^2 \omega_n^2}{EJ} \tag{25}$$

Table 1. End constraints for a single span Timoshenko beam with free ends.

| Constraint variable | Equivalent differential equations in $\phi_n(x)$ |
|---|---|
| $Q = 0$ | $\phi_n'''(x) + K_T \phi_n'(x) = 0$ |
| $M = 0$ | $\phi_n''(x) + K_M \phi_n(x) = 0$ |

Substituting the boundary condition requirements into the corresponding differential equations of the modal shapes, a system of simultaneous linear algebraic equations is obtained for each part of the frequency spectrum:

$$\mathbf{A}X = 0 \tag{26}$$

where, **A** is the coefficient matrix for the homogeneous system of linear algebraic equations, X is the unknown vector of the modal shape constants $C_i$:

$$X = \begin{bmatrix} C_1 \\ C_2 \\ C_3 \\ C_4 \end{bmatrix}, \quad 0 = \begin{bmatrix} 0 \\ 0 \\ 0 \\ 0 \end{bmatrix} \tag{27}$$

Table 2 indicates the coefficient matrix **A** for the homogeneous system of linear algebraic equations satisfying the free ends boundary conditions, for each part of the frequency spectrum.

For $X$ to be nonzero, the matrix **A** must be singular, i.e., its determinant must be zero, leading to the evaluation of the natural frequencies $\omega_n$ for the whole spectrum.

The proposed solution can be easily implemented within most programming languages, like Python (Van Rossum, 2015), allowing to evaluate the natural vibration frequencies and modes of a Timoshenko beam on Winkler foundation.

### 3.4. Modal equations

The solution of the partial differential equation Eq. (3) can be determined by solving the modal equations in the unknown modal coordinates $q_n(t)$ (Chopra, 1995):

$$M_n \ddot{q}_n(t) + C_n \dot{q}_n(t) + K_n q_n(t) = P_n(t) \tag{28}$$

where $M_n$, $C_n$, $K_n$ and $P_n(t)$ represent the generalized mass, damping, stiffness and force for the $n$th mode, respectively:

$$M_n = \int_0^L m_l [\phi_n(x)]^2 dx \tag{29}$$

$$C_n = 2\xi_n \omega_n M_n \tag{30}$$

$$K_n = \omega_n^2 M_n \tag{31}$$

$$P_n(t) = \int_0^L k_l \left( y_G(x,t) - \frac{EJ}{kGA_b} \frac{\partial^2 y_G(x,t)}{\partial x^2} + \frac{m_l r^2}{kGA_b} \frac{\partial^2 y_G(x,t)}{\partial t^2} \right) \phi_n(x) dx \tag{32}$$

where $\xi_n$ is the damping ratio for the $n$th mode, and $y_G(x,t) = U_g(x,t)$ is the ground displacement that determines the intensity of the soil-structure interaction (Figure 3).

The modal equations (Eq. (28)) can be solved numerically using the Newmark method for linear systems, permitting solution without iteration (Chopra, 1995).

Finally, the total beam displacement response $y(x,t)$ is calculated as a linear combination of the displacements of all vibration modes $y_n(x,t)$, using Eq. (4).

The main steps of the semi-analytical formulation are summarized in Table 3, outlining the progression from the governing differential equations to the total displacement response.

Table 2. Coefficient matrix **A** for the homogeneous system of linear algebraic equations satisfying the free ends boundary conditions, for different parts of the frequency spectrum.

| | **Frequency $\omega_n$** | **Coefficient matrix A** |
|---|---|---|
| 1. | $\omega_n < \tilde{\omega}_1$ | $\begin{bmatrix} a_{11} = -2\alpha\beta & a_{12} = \alpha^2 - \beta^2 + K_M & -a_{11} & a_{12} \\ a_{21} = -\beta^3 + 3\alpha^2\beta + \beta K_T & a_{22} = -\alpha^3 + 3\alpha\beta^2 - \alpha K_T & a_{21} & -a_{22} \\ e^{-\alpha L}\left(a_{12}\sin\beta L + a_{11}\cos\beta L\right) & e^{-\alpha L}\left(-a_{11}\sin\beta L + a_{12}\cos\beta L\right) & e^{\alpha L}\left(a_{12}\sin\beta L - a_{11}\cos\beta L\right) & e^{\alpha L}\left(-a_{11}\sin\beta L + a_{12}\cos\beta L\right) \\ e^{-\alpha L}\left(a_{22}\sin\beta L + a_{21}\cos\beta L\right) & e^{-\alpha L}\left(-a_{21}\sin\beta L + a_{22}\cos\beta L\right) & e^{\alpha L}\left(-a_{22}\sin\beta L + a_{21}\cos\beta L\right) & e^{\alpha L}\left(-a_{21}\sin\beta L - a_{22}\cos\beta L\right) \end{bmatrix}$ |
| 2. | $\omega_n = \tilde{\omega}_1$ | $\begin{bmatrix} 0 & b_{12} = -\lambda^2 + K_M & b_{13} = 2\lambda & 0 \\ b_{21} = -\lambda^3 + \lambda K_T & 0 & 0 & b_{24} = -3\lambda^2 + K_T \\ b_{12}\sin\lambda L & b_{12}\cos\lambda L & b_{12}L\sin\lambda L + b_{13}\cos\lambda L & b_{12}L\cos\lambda L - b_{13}\sin\lambda L \\ b_{21}\cos\lambda L & -b_{21}\sin\lambda L & b_{21}L\cos\lambda L + b_{24}\sin\lambda L & -b_{21}L\sin\lambda L + b_{24}\cos\lambda L \end{bmatrix}, \quad \lambda = \sqrt{\frac{B}{2A}}$ |
| 3. | $\tilde{\omega}_1 < \omega_n < \tilde{\omega}_2$ <br> $\vee$ <br> $\omega_n > \tilde{\omega}_3$ | $\begin{bmatrix} 0 & c_{12} = -\lambda_1^2 + K_M & 0 & c_{14} = -\lambda_2^2 + K_M \\ c_{21} = -\lambda_1^3 + \lambda_1 K_T & 0 & c_{23} = -\lambda_2^3 + \beta K_T & 0 \\ c_{12}\sin\lambda_1 L & c_{12}\cos\lambda_1 L & c_{14}\sin\lambda_2 L & c_{14}\cos\lambda_2 L \\ c_{21}\cos\lambda_1 L & -c_{21}\sin\lambda_1 L & c_{23}\cos\lambda_2 L & -c_{23}\sin\lambda_2 L \end{bmatrix}$ |
| 4. | $\omega_n = \tilde{\omega}_2 \vee$ <br> $\omega_n = \tilde{\omega}_3$ | $\begin{bmatrix} d_{11} = K_M & 0 & 0 & d_{14} = -\lambda_2^2 + K_M \\ 0 & d_{22} = K_T & d_{23} = -\lambda_2^3 + \lambda_2 K_T & 0 \\ d_{11} & d_{11}L & d_{14}\sin\lambda_2 L & d_{14}\cos\lambda_2 L \\ 0 & d_{22} & d_{23}\cos\lambda_2 L & -d_{23}\sin\lambda_2 L \end{bmatrix}, \quad \lambda_2 = \sqrt{\frac{B}{A}}$ |
| 5. | $\tilde{\omega}_2 < \omega_n < \tilde{\omega}_3$ | $\begin{bmatrix} 0 & e_{12} = \lambda_1^2 + K_M & 0 & e_{14} = -\lambda_2^2 + K_M \\ e_{21} = \lambda_1^3 + \lambda_1 K_T & 0 & e_{23} = -\lambda_2^3 + \lambda_2 K_T & 0 \\ e_{12}\sinh\lambda_1 L & e_{12}\cosh\lambda_1 L & e_{14}\sin\lambda_2 L & e_{14}\cos\lambda_2 L \\ e_{21}\cosh\lambda_1 L & e_{21}\sinh\lambda_1 L & e_{23}\cos\lambda_2 L & -e_{23}\sin\lambda_2 L \end{bmatrix}$ |

Table 3. Summary of the proposed semi-analytical framework

1. Obtain the fourth-order partial differential equation (PDE) of motion (Eq. 3), considering the equilibrium of a differential beam element under lateral excitation (Eq. 1) and the beam constitutive equations (Eq. 2).
2. Using the Fourier method of variable separation (Eq. 4), decompose the fourth-order partial differential equation (Eq. 3) in two ordinary partial differential equations, one governing the time function $q(t)$ (Eq. 5.1) and the other governing the spatial function $\phi(x)$ (Eq. 5.2).
3. Determine the roots $\lambda_{1/2}$ of the characteristic equation associated with the fourth-order PDE of the spatial function $\phi_n(\mathrm{x})$ (Eq. 5.2), as a function of the natural frequency $\omega_n$ (Eqs. 8 and 9).
4. Determine the transition frequencies $\tilde{\omega}_i$ (Eqs. 12-14) which define the range of the vibration spectrum where the two roots ($\lambda_{1/2}$) of the characteristic equation exhibit dinstict behaviour: (*i*) complex conjugates for $\omega_n < \tilde{\omega}_1$; (*ii*) purely imaginary roots for $\tilde{\omega}_1 < \omega_n < \tilde{\omega}_2 \vee \omega_n > \tilde{\omega}_3$; (*iii*) one real and one purely imaginary root for $\tilde{\omega}_2 < \omega_n < \tilde{\omega}_3$.
5. Determine the coefficient matrix **A** for the homogeneous system of linear algebraic equations satisfying the ends boundary conditions, for each part of the frequency spectrum (Table 2).
6. Determine the eigenfrequencies $\omega_n$, that render the determinant of the coefficient matrix **A** singular, along with the associated modal shapes $\phi_n(x)$, for each part of the frequency spectrum.
7. Determine the generalised mass $M_n$, damping $C_n$, stiffness $K_n$, and force $P_n(t)$ using Eqs. 29-32.
8. Solve the modal equations Eq. (28) numerically using the Newmark method for linear systems, permitting solution without iteration.
9. Calculate the total beam displacement response $y(x,t)$ as a linear combination of all modal contributions $y_n(x,t)$, using Eq. (4).

## 4. Dynamic analysis of a Timoshenko beam on Winkler Foundation

This section first introduces the buried-beam case study used to evaluate the accuracy of the proposed semi-analytical model, together with the finite-element formulation adopted for comparison. It then presents the free-vibration of the free–free Timoshenko beam on Winkler foundation, followed by the forced dynamic response under transient ground deformation (TGD). The combined analytical, numerical, and experimental evaluations provide deeper insight into the system's dynamic behaviour and demonstrate the validity of the proposed methodology.

### 4.1. Case study and finite element model

To investigate the dynamic response of a Timoshenko beam on elastic foundation, we consider the case study of a straight water pipeline with variable length and operational mass buried in medium dense sand. The considered API 120 X-42 grade steel pipeline, has a diameter $D$ = 107 cm (42 in.) and wall thickness $t$ = 8 mm (5/16 in.) and is assumed buried in medium dense sand at a depth $H$ = 1.65 m, which is measured from the soil surface to the pipe axis. The soil–structure system parameters used in this study are summarized in Table 4.

The X42 steel grade pipe material is defined as elastic, with Young modulus $E$ = 210 GPa, and Poisson's ratio $\nu = 0.3$. The soil-pipeline interaction in the axial and transverse horizontal direction is modelled with uniaxial springs, calculated according to the ALA-ASCE guidelines (ASCE 1984; ALA 2001), considering medium dense sand soil conditions, with friction angle $\phi = 30°$ and density $\gamma$ = 16.9 kN/m$^3$ (Table 4). Within these guidelines, the force-displacement relationship of the soil springs, or p-y curves, is defined in terms of elasto-plastic functions, obtained by fitting the hyperbolic equation proposed by Trautmann and O'Rourke (1983). The elasto-plastic model is fully characterized by the maximum resistance $p_u$ and the deformation $y_u$ at yield, and the 1984 ASCE Guidelines suggest using as the effective stiffness twice the ratio of ultimate resistance to the yield deformation, $k_l = 2p_u/y_u$ (ASCE 1984; O'Rourke and Liu, 2012).

In this case, the calculated soil friction reaction per unit length of pipe is $f_r$ = 41.22 kN/m, while the maximum lateral soil reaction is $p_u$ = 1365.22 kN/m, achieved at a relative soil-pipe displacement $u_0$ = 0.15 cm and $y_u$ = 10.91 cm, respectively. The resulting soil stiffness in the axial and transverse direction is $k_a$ = 56.0 MPa and $k_l$ = 2.5 MPa, respectively (Table 4).

Table 4. Soil-structure system parameters.

| | | |
|---|---|---|
| Pipeline diameter, $D$ | 1.067 | m |
| Pipeline wall thickness, $t$ | 0.0079 | m |
| Pipeline burial depth, $H$ | 1.65 | m |
| Soil unit weight, $\gamma$ | 16.9 | kN/m$^3$ |
| Soil friction angle, $\phi$ | 30 | ° |
| Soil friction reaction per unit pipe length, $f_r$ | 41.22 | kN/m |
| Axial soil-spring maximum elastic deformation, $u_0$ | 0.15 | cm |
| Linear axial soil-spring stiffness, $k_a = 2f_r/u_0$ | 56.044 | MPa |
| Lateral soil reaction per unit pipe length, $p_u$ | 1365.22 | kN/m |
| Lateral soil-spring maximum elastic deformation, $y_u$ | 10.91 | cm |
| Lateral soil-spring stiffness, $k_l = 2p_u/y_u$ | 2.503 | MPa |
| Steel elastic modulus, $E$ | 210 | GPa |
| Steel Poisson's ratio, $\nu$ | 0.3 | |
| Elastic shear modulus, $G = 0.5E/(1+\nu)$ | 80.77 | GPa |
| Beam shear coefficient, $k$ | 0.53 | |
| Pipeline cross-sectional area, $A_b$ | 0.026 | m$^2$ |
| Beam second moment of inertia, $J$ | 0.0037 | m$^4$ |
| Radius of gyration of the beam cross section, $r = \sqrt{J/A_b}$ | 0.374 | m |
| Steel density, $\rho_s$ | 7860.35 | kg/m$^3$ |
| Water density, $\rho_w$ | 1000.00 | kg/m$^3$ |
| Linear mass of the unfilled pipeline, $m_l$ | 207.56 | kg/m |
| Linear mass of the filled water pipeline, $m_{lw}$ | 1074.99 | kg/m |

The numerical dynamic modal analysis of the buried pipeline subjected to ground vibration is performed using the finite element software ABAQUS/Standard (Simulia, 2020). The buried pipeline was modeled using PIPE21H beam elements, while the soil-structure interaction was represented with the spring-like pipe-soil interaction elements PSI24, as schematically illustrated in Figure 5.

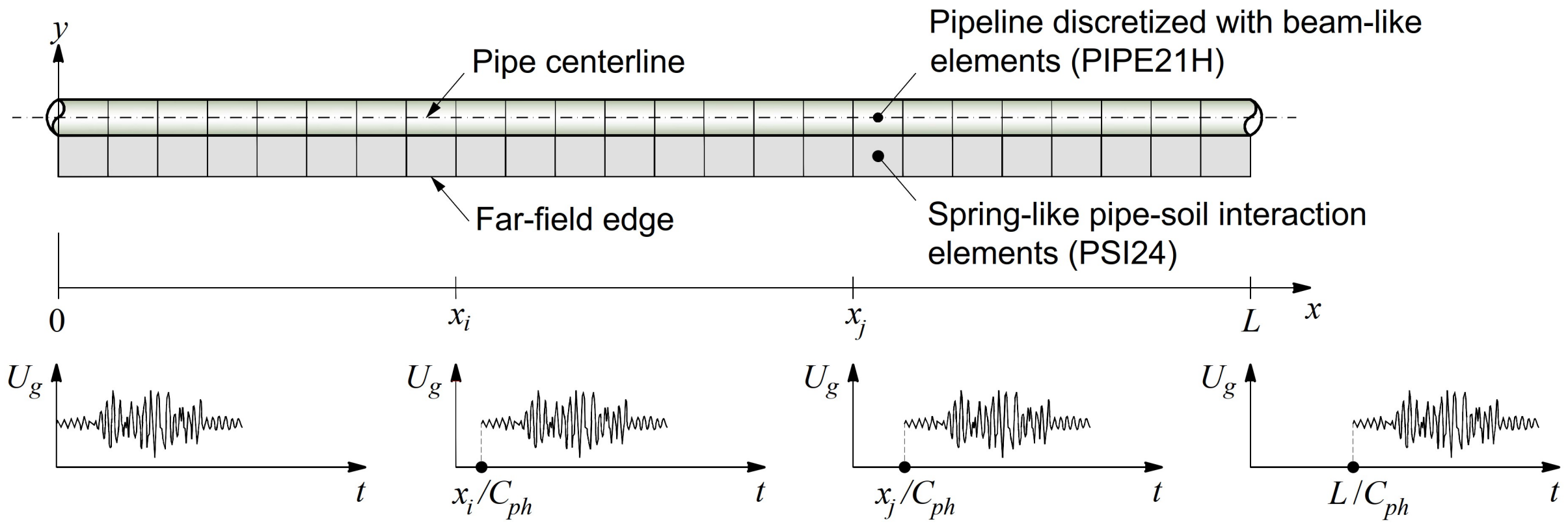

Figure 5. Finite element model of the pipeline-soil system subjected to seismic induced ground displacement time histories $U_g(x,t)$, with a time shift at axial distance $x_i$ equal to $x_i/C_{ph}$, inversely proportional to the apparent wave velocity $C_{ph}$.

One edge of the soil element shares nodes with the overlying pipe element while the nodes on the other edge are assigned the base ground motion with a time shift proportional to the distance along the system axis (Figure 5). A 1.0 m mesh size was adopted for both pipe and soil elements, based on a mesh sensitivity study ensuring solution accuracy and computational efficiency.
The numerical dynamic modal analysis is conducted in two consecutive steps. First, the eigenvalue extraction of the first 1000 modes is performed, calculating the corresponding natural frequencies and the mode shapes of the system. Second, in the modal dynamic step, the modal amplitudes are integrated through time, and the dynamic response is obtained by modal superposition. The proposed modal dynamic analysis procedure has the advantage of being computationally more efficient than direct-integration methods, while providing useful insight into the system dynamic behaviour.

### 4.2. Natural vibration frequencies and modes of Timoshenko beam on Winkler Foundation

The evaluated natural frequencies $\omega_n$ using the proposed analytical solution and the finite element modal analysis for the first 20 vibration modes are indicated in Table 5. The percent difference between the two approaches does not exceed 0.2%, demonstrating the validity of the proposed analytical solution.

Table 5. Computed natural frequencies $\omega_n$ (rad/s) and mode numbers $n$ of the buried pipelines of different lengths $L$, represented as Timoshenko beams on Winkler foundation, according to the semi-analytical (SA) and the finite element model (FEM).

| | $L$ = 100 m | | | $L$ = 1000 m | | | $L$ = 10000 m | | |
|---|---|---|---|---|---|---|---|---|---|
| n | SA (rad/s) | FEM (rad/s) | % difference | SA (rad/s) | FEM (rad/s) | % difference | SA (rad/s) | FEM (rad/s) | % difference |
| 1 | 109.804 | 109.805 | 0.001 | 109.808 | 109.811 | 0.003 | 109.810 | 109.811 | 0.002 |
| 2 | **109.813** | **109.811** | 0.002 | 109.810 | 109.811 | 0.001 | 109.812 | 109.811 | 0.001 |
| 3 | 109.865 | 109.861 | 0.003 | 109.812 | 109.811 | 0.001 | 109.812 | 109.811 | 0.001 |
| 4 | 110.385 | 110.370 | 0.013 | 109.812 | 109.811 | 0.001 | 109.812 | 109.811 | 0.001 |
| 5 | 112.134 | 112.111 | 0.020 | 109.813 | 109.811 | 0.001 | 109.812 | 109.811 | 0.001 |
| 6 | 116.147 | 116.107 | 0.035 | 109.813 | 109.811 | 0.001 | 109.812 | 109.811 | 0.001 |
| 7 | 123.537 | 123.471 | 0.053 | 109.813 | 109.811 | 0.002 | 109.812 | 109.811 | 0.001 |
| 8 | 135.210 | 135.101 | 0.081 | **109.813** | **109.811** | 0.002 | 109.812 | 109.811 | 0.001 |
| 9 | 151.668 | 151.506 | 0.107 | 109.813 | 109.811 | 0.002 | 109.812 | 109.811 | 0.001 |
| 10 | 172.987 | 172.769 | 0.126 | 109.815 | 109.818 | 0.002 | 109.812 | 109.811 | 0.001 |
| 11 | 198.948 | 198.674 | 0.138 | 109.819 | 109.818 | 0.001 | 109.812 | 109.811 | 0.001 |
| 12 | 229.197 | 228.878 | 0.139 | 109.824 | 109.824 | 0.000 | 109.812 | 109.811 | 0.001 |
| 13 | 263.347 | 263.002 | 0.131 | 109.831 | 109.830 | 0.001 | 109.812 | 109.811 | 0.001 |
| 14 | 301.034 | 300.688 | 0.115 | 109.840 | 109.843 | 0.002 | 109.812 | 109.811 | 0.001 |
| 15 | 341.928 | 341.617 | 0.091 | 109.853 | 109.855 | 0.002 | 109.812 | 109.811 | 0.001 |
| 16 | 385.737 | 385.505 | 0.060 | 109.869 | 109.868 | 0.001 | 109.812 | 109.811 | 0.001 |
| 17 | 432.201 | 432.088 | 0.026 | 109.889 | 109.887 | 0.002 | 109.812 | 109.811 | 0.001 |
| 18 | 481.086 | 481.147 | 0.013 | 109.913 | 109.912 | 0.001 | 109.812 | 109.811 | 0.001 |
| 19 | 532.181 | 532.469 | 0.054 | 109.943 | 109.943 | 0.000 | 109.812 | 109.811 | 0.001 |
| 20 | 585.293 | 585.863 | 0.097 | 109.979 | 109.981 | 0.002 | 109.812 | 109.811 | 0.001 |

Figure 6 illustrates the full frequency spectrum for the first $N$ = 100 modes, highlighting the number of modes for each part of the spectrum, considering the non operating and operating pipeline filled with water. The former presents greater vibration frequencies than the latter, because of the lower mass of the unfilled pipeline. It is important to note that the first part of the spectrum contains only one vibration mode for the shorter beam ($L$ = 100 m) and two distinct modes for the two longer beams. The remaining parts of the spectrum contain a density of vibration modes that increases with the pipeline length $L$.
Particularly, the second part of the frequency spectrum contains only the rigid-body vibration mode for the shorter pipeline, $L$ = 100 m ($N_2$ = 1), and an increasing number of modes for the 1000 m ($N_2$ = 6) and the 10000 m ($N_2$ = 66) long pipelines. Clearly, the vibration frequencies become the closer to each other, the longer the system's length (Table 5). This is consistent with previous research publications, noting that, in longer pipelines, the natural frequencies of all lateral modes become closer to each other, requiring a significantly greater number of modes to calculate dynamic response (Hindy and Novak, 1980; Zerva et al., 1988; Hosseini and Roudsari, 2015).

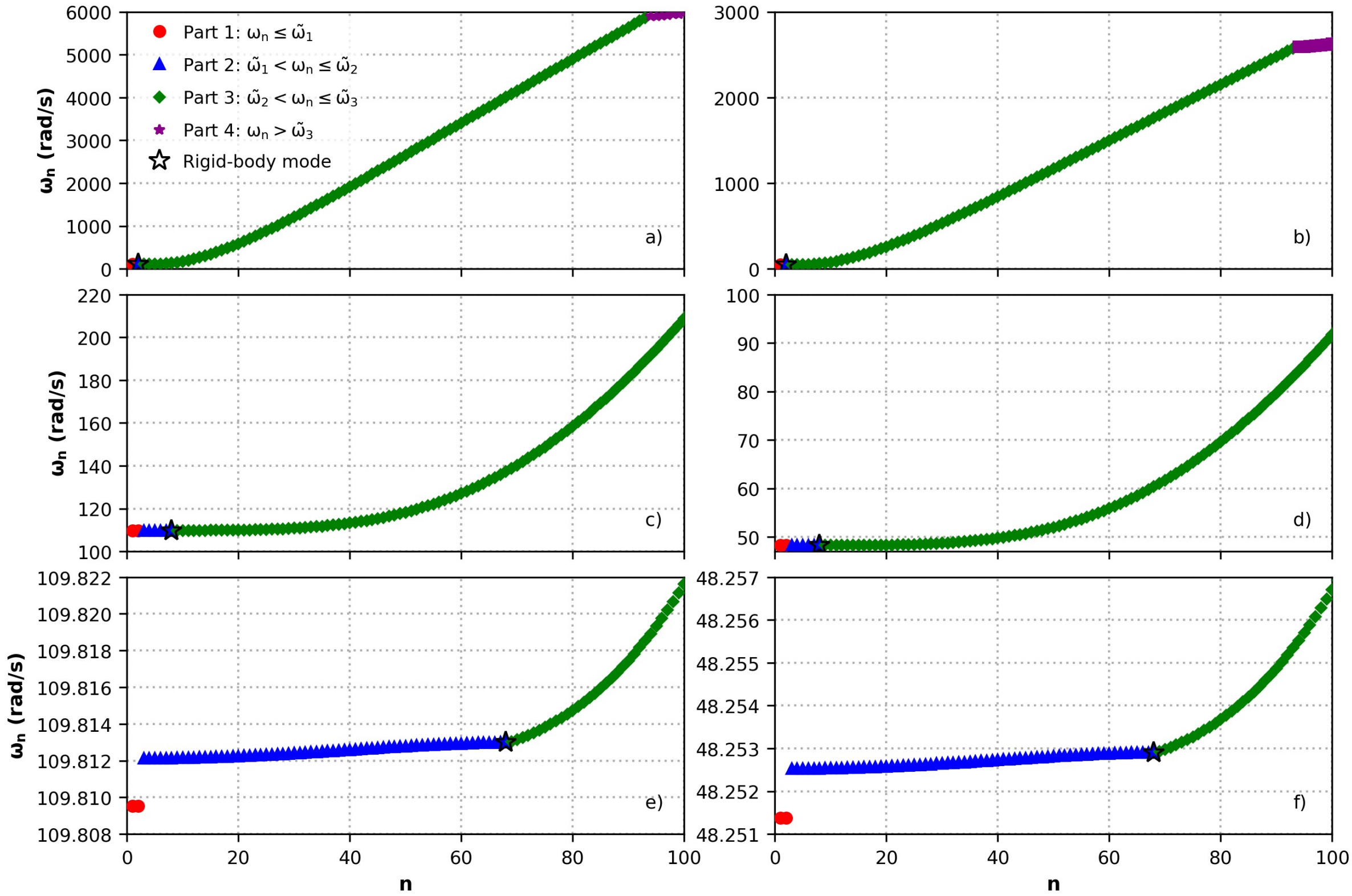

Figure 6. Frequency spectrum representing the first 100 vibration frequencies $\omega_n$ as a function of the wave number $n$ for different beam lengths, and operating conditions: $L$ = 100 m, a) unfilled and b) filled with water; $L$ = 1000 m: c) unfilled and d) filled with water; $L$ = 10000 m: e) unfilled and f) filled with water.

The cut-off frequencies $\tilde{\omega}_1$ (109.812 rad/s) and $\tilde{\omega}_2$ (109.813 rad/s), delimiting the second part of the frequency spectrum, are very close between each other. The second cut-off frequency $\tilde{\omega}_2$, corresponding to the rigid body mode (Eq. (13)), is equal to the square root of the lateral soil stiffness-to-mass ratio ($k_l/m_l$), characterizing the system dynamic behaviour. If the forcing frequency generated by the vibration source is close to the first two cut-off frequencies, the system will undergo dynamic amplification.

Figure 7 shows the variation of the second cut-off frequency $\tilde{\omega}_2$ as the soil stiffness is reduced from 2.5 MPa to lower values for the cases of unfilled pipeline and water filled pipeline. As the lateral soil stiffness decreases due to poor soil compaction, the cut-off frequency reduces. Moreover, the increased effective mass of the buried pipe carrying water tends to further reduce its natural frequencies. In this pipeline case, the system is more vulnerable to dynamic amplification effects caused by earthquake vibration frequencies of a few Hz.
Particularly, at resonance, when the forcing frequency $\omega$ is close to the system natural frequency $\omega_n$, the predominant deformation response of the structure is determined by the corresponding modal shape $\phi_n(x)$.

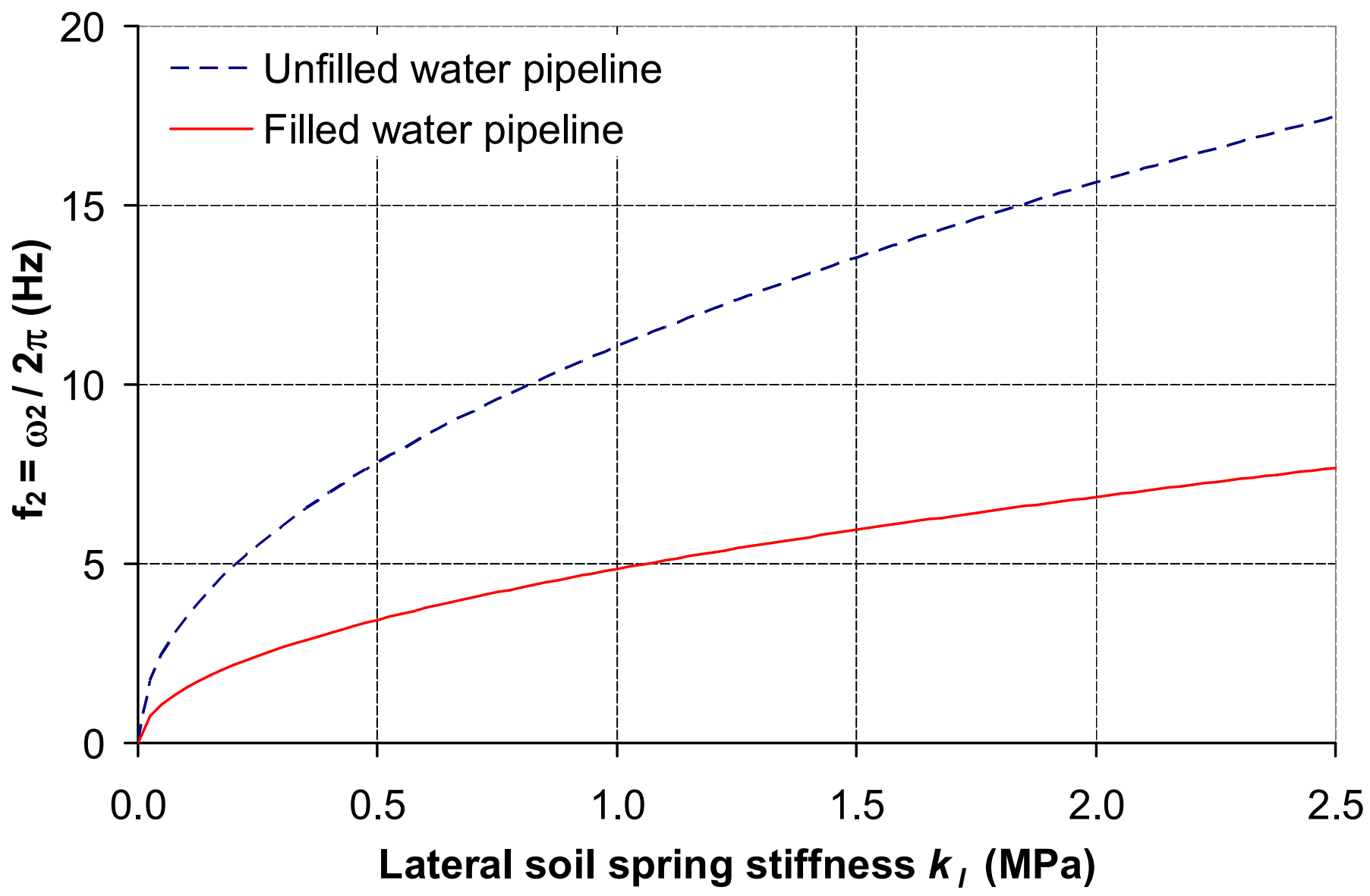

Figure 7. Cut-off frequency $\tilde{f}_2 = \tilde{\omega}_2 / 2\pi$ as a function of the lateral stiffness $k_l$, for different effective mass $m_l$ associated with each pipeline operating condition.

The evaluated modal shapes for the first, second, and third part of the frequency spectrum, are shown in Figures 8, 9 and 10, respectively, considering the proposed analytical and finite element modal analysis.
Comparison between the semi-analytical and numerical results shows excellent agreement between the two approaches, further confirming the validity of the proposed semi-analytical model.
For the shorter buried beam ($L$ = 100 m), the first shape function $\phi_1(x)$ becomes similar to a rigid body rotation mode due to the minimal curvature (Figure 8a). This is consistent with results reported in other publications on free vibration of short beams on elastic foundation as well (Wang, 1978; Wu, 2019). Conversely, the longer beams present two distinct curved shape functions that diverge significantly from a rigid body rotation, according to Eq. (18), as illustrated in Figure 8b and 8c for $L$ = 1000 m and $L$ = 10000 m, respectively.
Within the second part of the frequency spectrum ($\tilde{\omega}_1 < \omega_n \leq \tilde{\omega}_2$), the shorter beam manifests only the rigid body translational mode (Figure 9a). On the other hand, the 1000 m and 10000 m long beam exhibit additional 5 and 66 harmonic shape functions, respectively.
In the third part of the spectrum ($\tilde{\omega}_2 < \omega_n \leq \tilde{\omega}_3$), the third shape function $\phi_3(x)$ of the shorter buried beam ($L$ = 100 m) corresponds to a global bending mode, whereas the subsequent modal shapes exhibit increasing inflection points along its axis (Figure 10a). Conversely, for the longer beams, this part of the spectrum includes mode numbers greater than 8 and 69 for the 1000 m (Figure 10b) and 10000 m (Figure 10c) long system respectively, resulting in increased number of inflection points and shorter average wave lengths (< 285 m).
In the first and third part of the frequency spectrum, the modal shape functions assume their maximum values at the beam's end, as particularly evident for the 100 m long system (Figure 8a and 10a). This leads to a greater dynamic amplification herein at resonance, when the forcing frequency is close to the system's natural frequency in that spectrum range.

Because of the higher frequency and greater curvature, a greater energy is required to excite higher modes, that tend to contribute less to the total vibration amplitude, particularly for shorter beams, characterized by a lower modal density. Conversely, longer beams exhibit a higher modal density, with the natural frequencies closer to each other, which causes excitation of a greater number of modes contributing to the dynamic amplification response.

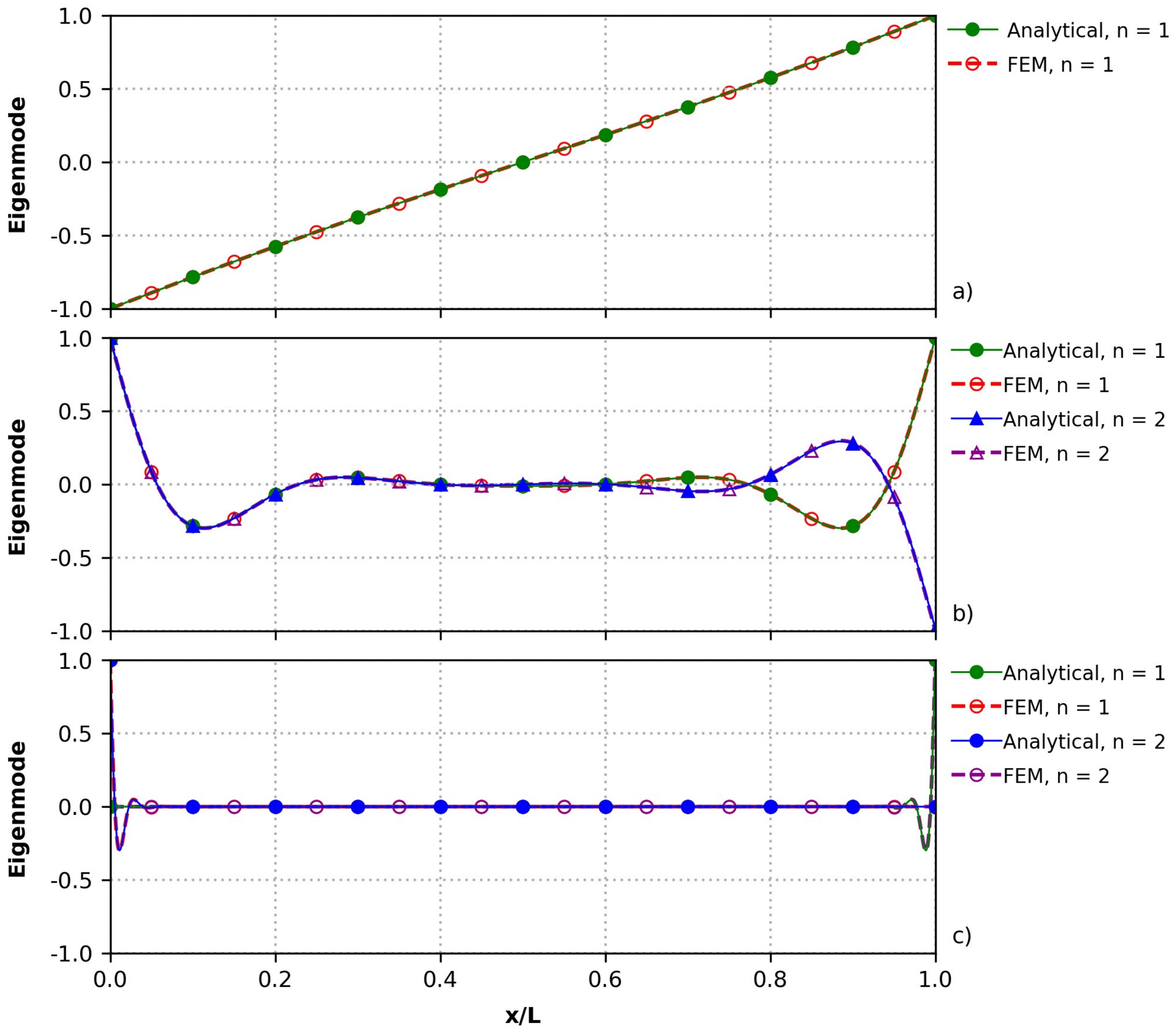

Figure 8. Natural vibration modes for the first part of the frequency spectrum ($\omega_n \leq \tilde{\omega}_1$) for various pipeline lengths: a) $L$ = 100 m; b) $L$ = 1000 m; c) $L$ = 10000 m.

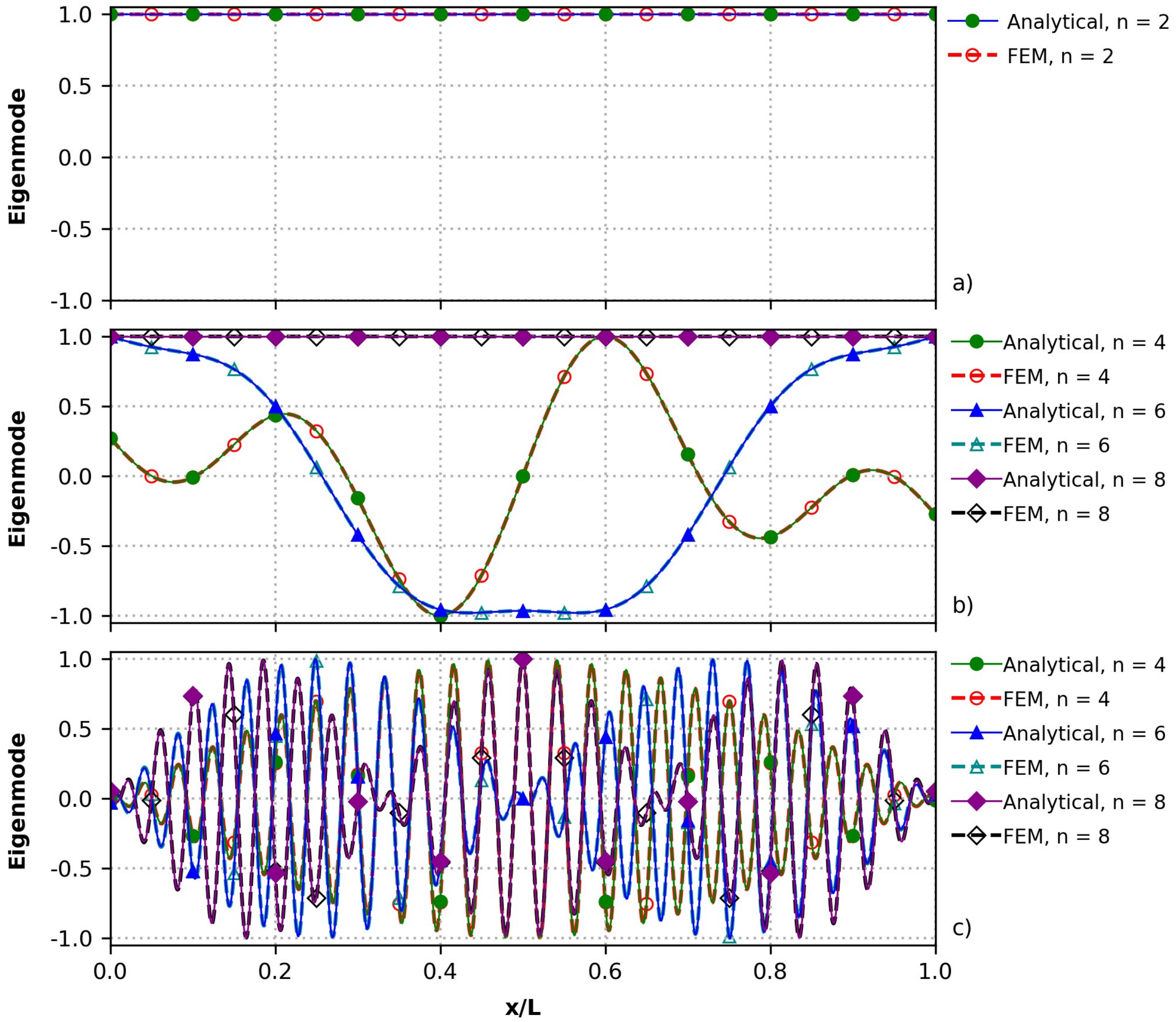

Figure 9. Lowest vibration modes for the second part of the frequency spectrum ($\tilde{\omega}_1 < \omega_n \leq \tilde{\omega}_2$) corresponding to even mode numbers $n$, for various pipeline lengths: a) $L$ = 100 m; b) $L$ = 1000 m; c) $L$ = 10000 m.

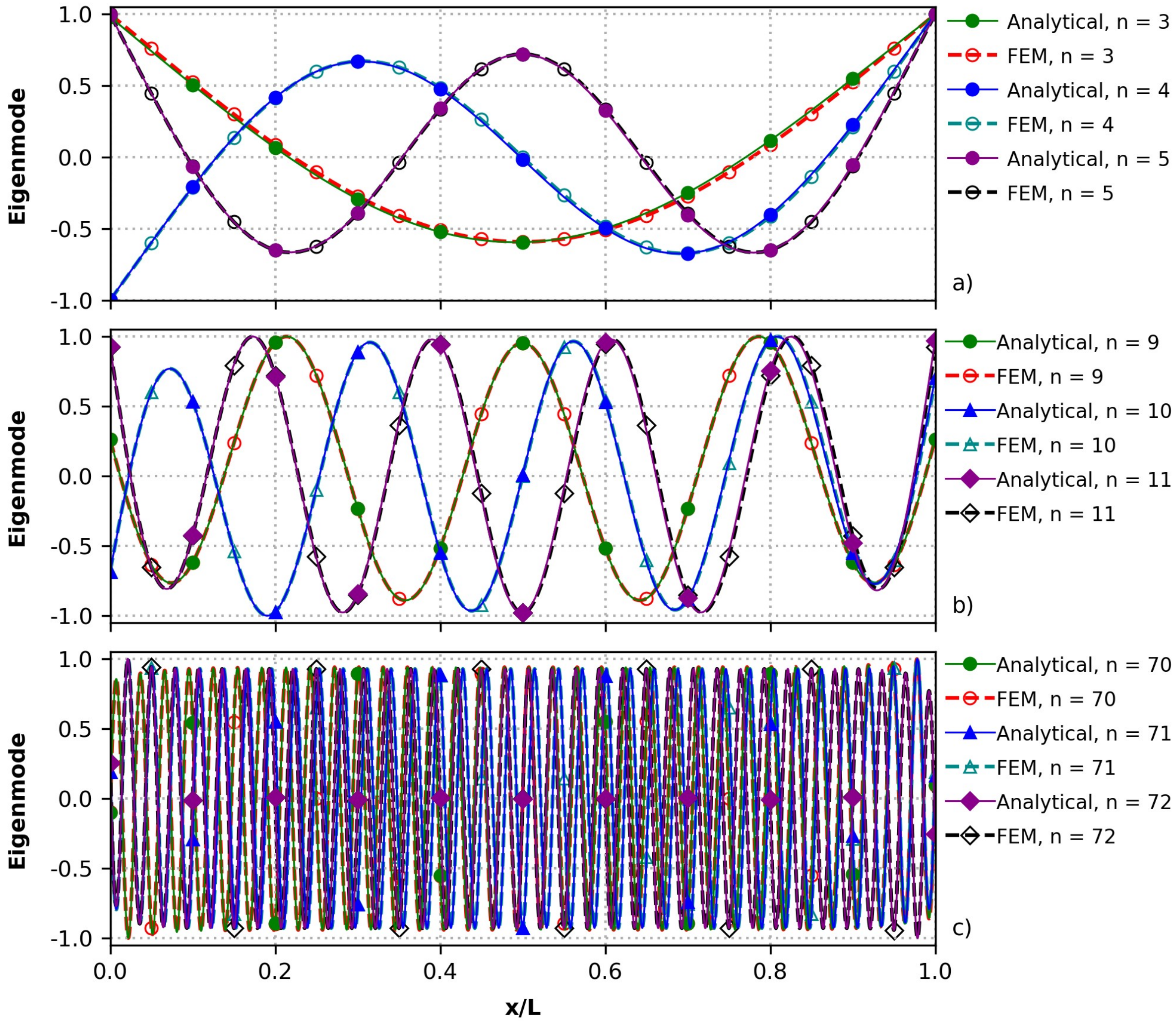

Figure 10. Lowest vibration modes for the third part of the spectrum ($\tilde{\omega}_2 < \omega_n \leq \tilde{\omega}_3$) for various pipeline lengths: a) $L = 100$ m; b) $L = 1000$ m; c) $L = 10000$ m.

### 4.3. Experimental validation

This section evaluates the accuracy of the proposed semi-analytical formulation through two complementary experimental benchmarks. First, the free-vibration behavior of a free-free Timoshenko beam is examined using high-precision laboratory measurements reported by Monsivais et al. (2016), representing the limiting case of vanishing foundation stiffness. Second, the forced-vibration response of a buried pipeline subjected to transverse earthquake-induced ground excitation is assessed using multi-point shaking-table tests by Han et al. (2024). Together, these experimental datasets validate the model across a broad range of system properties and dynamic conditions.

#### 4.3.1. Experimental validation for the free vibration response of a free-free Timoshenko beam

This section compares the free vibration response of a free-free Timoshenko beam evaluated using the proposed analytical method and the high precision experimental results reported in Monsivais et al. (2016), representing the limit case of vanishing Winkler foundation ($k_l \to 0$). The authors confirmed experimentally the Timoshenko beam theory prediction of a critical transition frequency $\tilde{f}_c$, dividing the frequency spectrum in two distinct parts of different oscillatory characteristics. The tested beam had a length $L$ = 0.500 m, height $a$ = 0.0252 m, width $b$ = 0.0504 m, material density $\rho$ = 2700 kg/m$^3$, Young modulus $E$ = 67.421 GPa, shear modulus $G$ = 26.919 GPa and a Timoshenko shear coefficient $k$ = 0.85.

Figure 11 compares the first 12 natural frequencies $f_n$ measured by Monsivais et al. (2016) with the predictions of the present analytical method (excluding rigid-body modes) for the first ($f_n < \tilde{f}_c$) and the second part of the frequency spectrum ($f_n > \tilde{f}_c$), showing excellent agreement. As shown in Figure 12, the modal shapes evaluated using the present analytical model agree very well with the high resolution amplitude measurements by Monsivais et al. (2016), further confirming the validity of the proposed procedure.

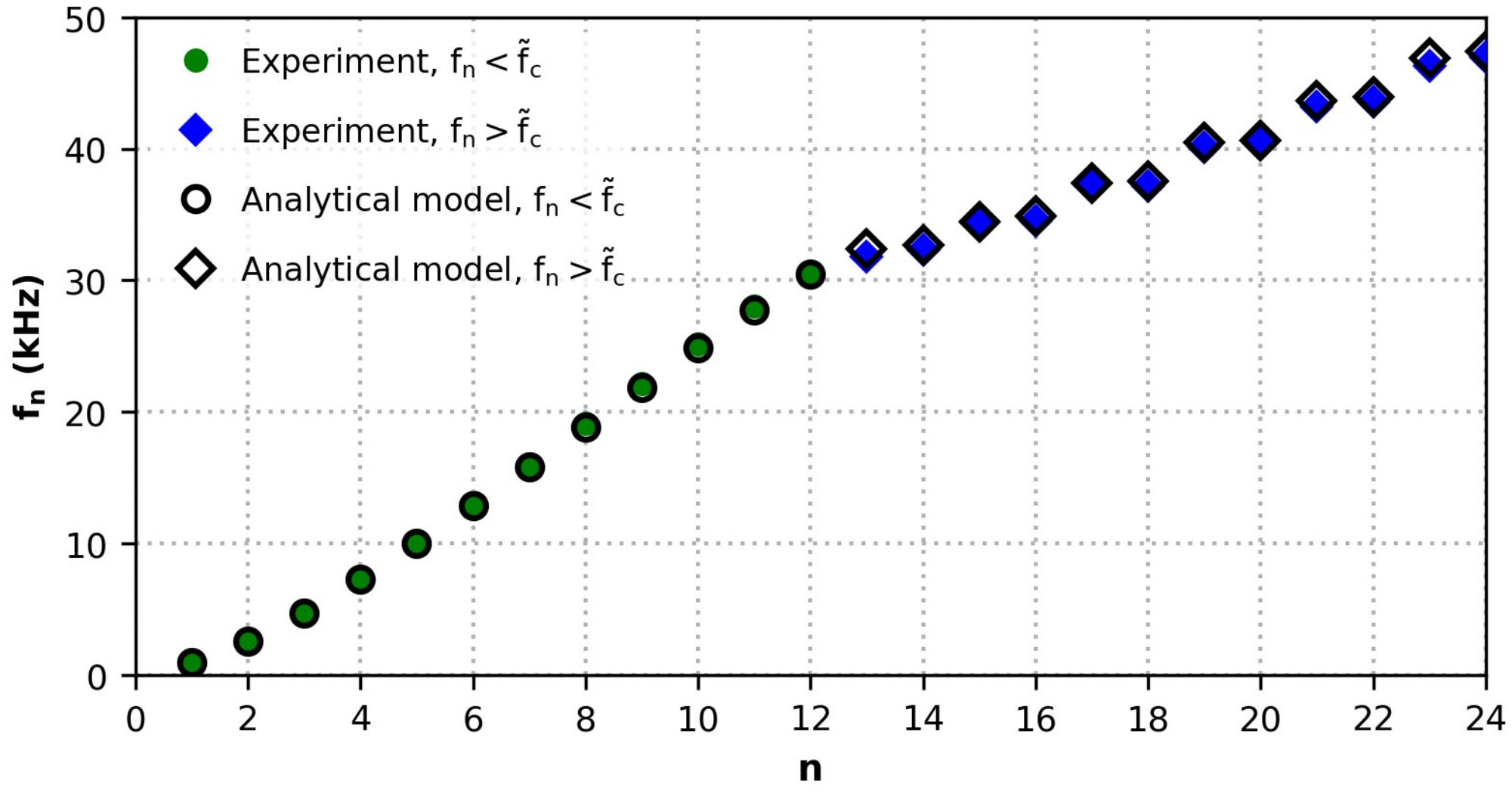

Figure 11. Frequency spectrum representing the first 24 vibration frequencies $f_n$ as a function of the mode number $n$ (excluding rigid-body modes) for the bare free-free Timoshenko beam; comparison between the predictions of the present analytical method and the free-vibration tests on Timoshenko beams by Monsivais et al. (2016) for the first ($f_n < \tilde{f}_c$) and the second part of the frequency spectrum ($f_n > \tilde{f}_c$).

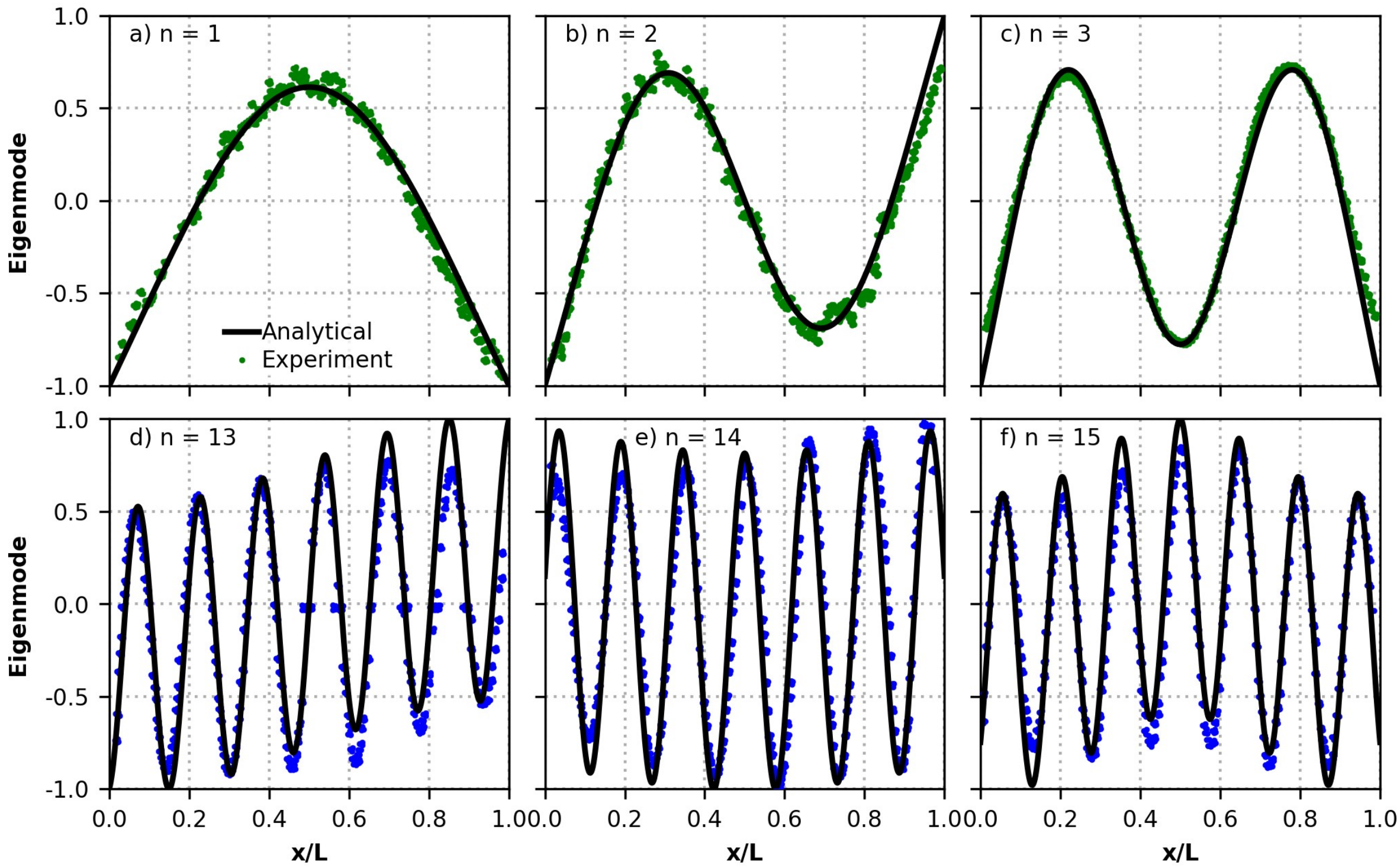

Figure 12. Lowest three natural vibration modes for each part of the frequency spectrum, according to the proposed analytical method (continuous line) and the experimental amplitude measurements by Monsivais et al. (2016), given by dots.

#### 4.3.2. Experimental validation for the forced vibration response of a buried Timoshenko beam under transverse ground excitation

The present methodology is further validated considering recent shaking table tests results, employing recorded earthquake ground motions in the transverse direction to investigate the seismic response of a buried pipeline (Han et al., 2024). The tested pipe had a length $L$ = 6.0 m, diameter $D$ = 0.15 m, thickness $t$ = 0.002 m, material density $\rho$ = 2100 kg/m$^3$ and Young modulus $E$ = 3.2 GPa. It is embedded in medium dense sand of density $\gamma$ = 1740 kg/m$^3$, at a burial depth $H$ = 0.275 m.

Figure 13 compares the acceleration time history response of the buried pipeline in the transverse direction under the 0.25 g EL Centro uniform excitation, obtained using the present semi-analytical method and the shaking table test measurements by Han et al. (2024).

The analytical and experimental responses show very good agreement throughout the time history, with peak acceleration differences below 3.3%, confirming the accuracy of the proposed semi-analytical method.

Together, the two validation cases confirm that the proposed semi-analytical framework accurately predicts both the free and forced dynamic behavior of buried Timoshenko beams across a wide frequency range.

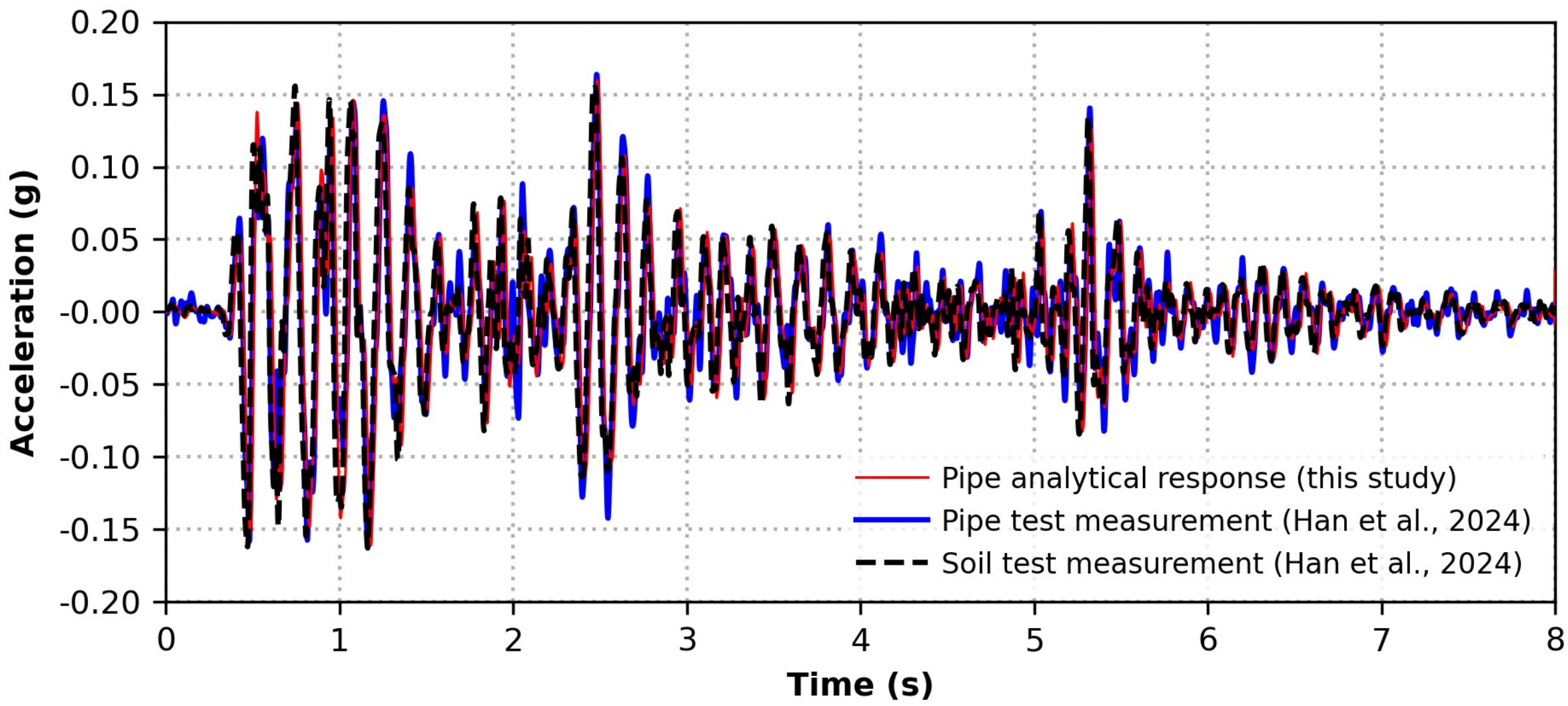

Figure 13. Acceleration response of the buried pipeline in the transverse direction under the 0.25 g EL Centro uniform excitation; comparison between the predictions of the present methodology and the shaking table test measurements by Han et al. (2024).

### 4.4. Forced dynamic response of the Timoshenko beam on Winkler foundation

In this case study, the ground motion $U_g(x,t)$ is modelled as a sinusoidal shear wave with amplitude $U_{g,max}$, propagating with an angular frequency $\omega$ and apparent shear wave velocity $C_{ph}$ = 2000m/s. This value is within with the range of horizontal propagation velocity values reported in O'Rourke and El Hmadi (1988), Kiyomiya (1995). The total duration of the ground motion is considered equal to $t_g$ = 30 s.

To take into account the wave propagation effect, the imposed ground motion is applied at the free ends of the soil springs at a harmonic function of the distance $x$ along the pipe axis, with a time lag equal to $x/C_{ph}$:

$$U_g(x,t) = \begin{cases} U_{g,\max} \sin \omega \left( t - \dfrac{x}{C_{ph}} \right), & x / C_{ph} \leq t \leq x / C_{ph} + t_g \\ 0, \quad t < x / C_{ph} \ \vee \ t > x / C_{ph} + t_g \end{cases} \tag{33}$$

Although a harmonic travelling wave is used here to illustrate frequency-dependent behaviour, the underlying semi-analytical solution is fully general and applies to any time-history ground motion input, including recorded seismic data.

The external force in the right side of the differential equation of motion (Eq. (3)) is directly proportional to the ground displacement $y_G(x,t) = U_g(x,t)$ that determines the intensity of the soil-pipeline interaction.

The forced dynamic response of the Timoshenko beam on Winkler foundation is evaluated by solving the modal equations in the unknown modal coordinates numerically and then calculating the total displacement as a linear combination of all modes (Table 3). In the semi-analytical approach, 10, 90, and 230 modes were required for the 100 m, 1000 m, and 10000 m long system, respectively. The proposed algorithm can be easily implemented within most programming languages, like Python (Van Rossum, 2015), allowing to evaluate the dynamic behaviour of a Timoshenko beam on Winkler foundation under transient ground deformation.

Figure 14 shows the relative pipeline displacement, $U_p/U_{g,max}$, along its axis at $t$ = 5 s, obtained using the semi-analytical, finite element modal and implicit dynamic analysis, assuming a ground motion vibration frequency $f$ = 0.5. Excellent agreement between the analysis approaches is apparent.

The total pipeline vibration is the sum of two harmonic curves, characterized by the ground motion forcing frequency $f = 0.5$ Hz, and the fundamental frequency of the buried pipeline ($f_1 = \omega_1/2\pi = 17.5$ Hz). The former is smaller than the latter, leading to a quasi-static system response. The pipeline generally follows the ground movement, with the ratio of the maximum pipe displacement to the maximum ground displacement, $U_{p,max} / U_{g,max}$, equal to 1.09, 1.12, and 1.11 for the 100 m, 1000 m, and 10000 m long system, respectively. This ratio remains close to unity in all cases and is slightly lower for the shortest pipeline, whose length is shorter than the wavelength associated with the system's fundamental frequency ($\lambda_1 = C_{ph}/f_1 = 114$ m).

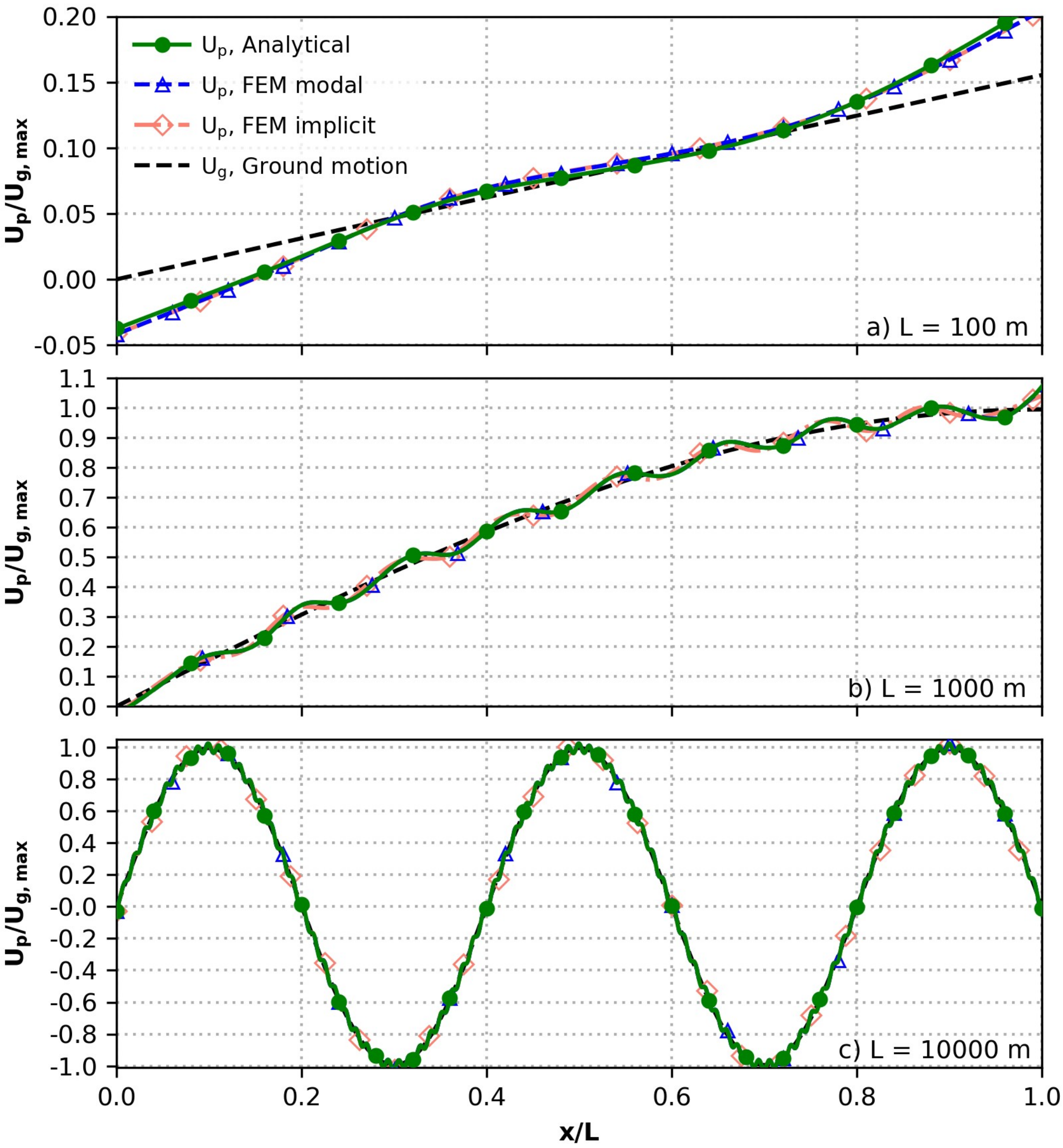

.

Figure 14. Relative pipeline displacement $U_p/U_{g,max}$ along its axis ($x/L$) at time $t$ = 5 s according to the semi-analytical, finite element based modal analysis and implicit dynamic analysis, considering an input sinusoidal ground displacement with frequency $f = 0.5$ Hz, and apparent wave propagation velocity $C_{ph}$ = 2000 m/s, for various pipeline lengths: a) $L$ = 100 m; b) $L$ = 1000 m; c) $L$ = 10000 m.

Figure 15 illustrates the variation of the lateral displacement ratio $U_{p,max} / U_{g,max}$ as a function of time ($t$) and space ($x/L$), for different lengths of the unfilled water pipeline, considering an input sinusoidal ground displacement with frequency $f = 0.5$ Hz. The pipeline follows the harmonic ground displacement quasi-statically, with the ratio $U_{p,max} / U_{g,max}$ approaching unity for longer pipeline lengths, throughout the duration of the ground motion. However, as discussed earlier, a smaller soil stiffness associated with poorly compacted backfill soil, as well the greater inertia of the filled water pipe, tend to reduce the system natural frequencies, making it more vulnerable to seismic dynamic amplification.

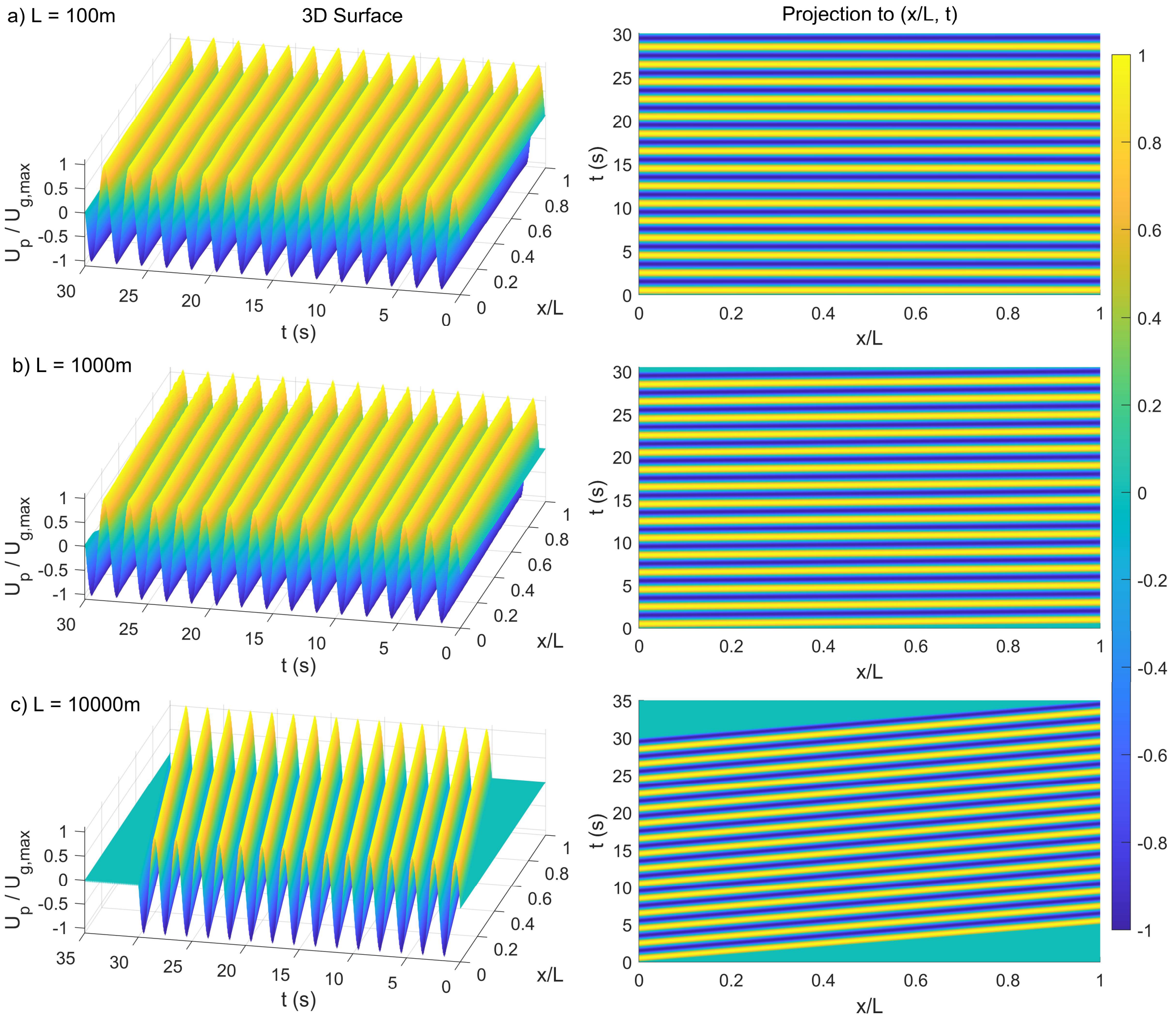

Figure 15. Variation of the relative pipeline displacement $U_p/U_{g,max}$ along its axis ($x/L$) during the seismic wave propagation with $f = 0.5$, for various pipeline lengths: a) $L = 100$ m; b) $L = 1000$ m; c) $L = 10000$ m.

Figure 16 illustrates the variation of the lateral displacement ratio $U_{p,max} / U_{g,max}$ as a function of time ($t$) and space ($x/L$), for different lengths ($L$) of the filled water pipeline and reduced soil stiffness ($k_l = 11$ kN/m$^2$), with fundamental frequency close to that of the ground motion ($f = 0.5$ Hz). Conversely to the previous case of unfilled water pipeline embedded in poorly compacted backfill, the displacement ratio $U_{p,max} / U_{g,max}$ increases with time, reaching values of 47.2, 47.6, and 49.0, for the 100m, 1000m and 10000m long system, respectively. The system has reached resonance, with significant amplification of the soil-structure interaction.

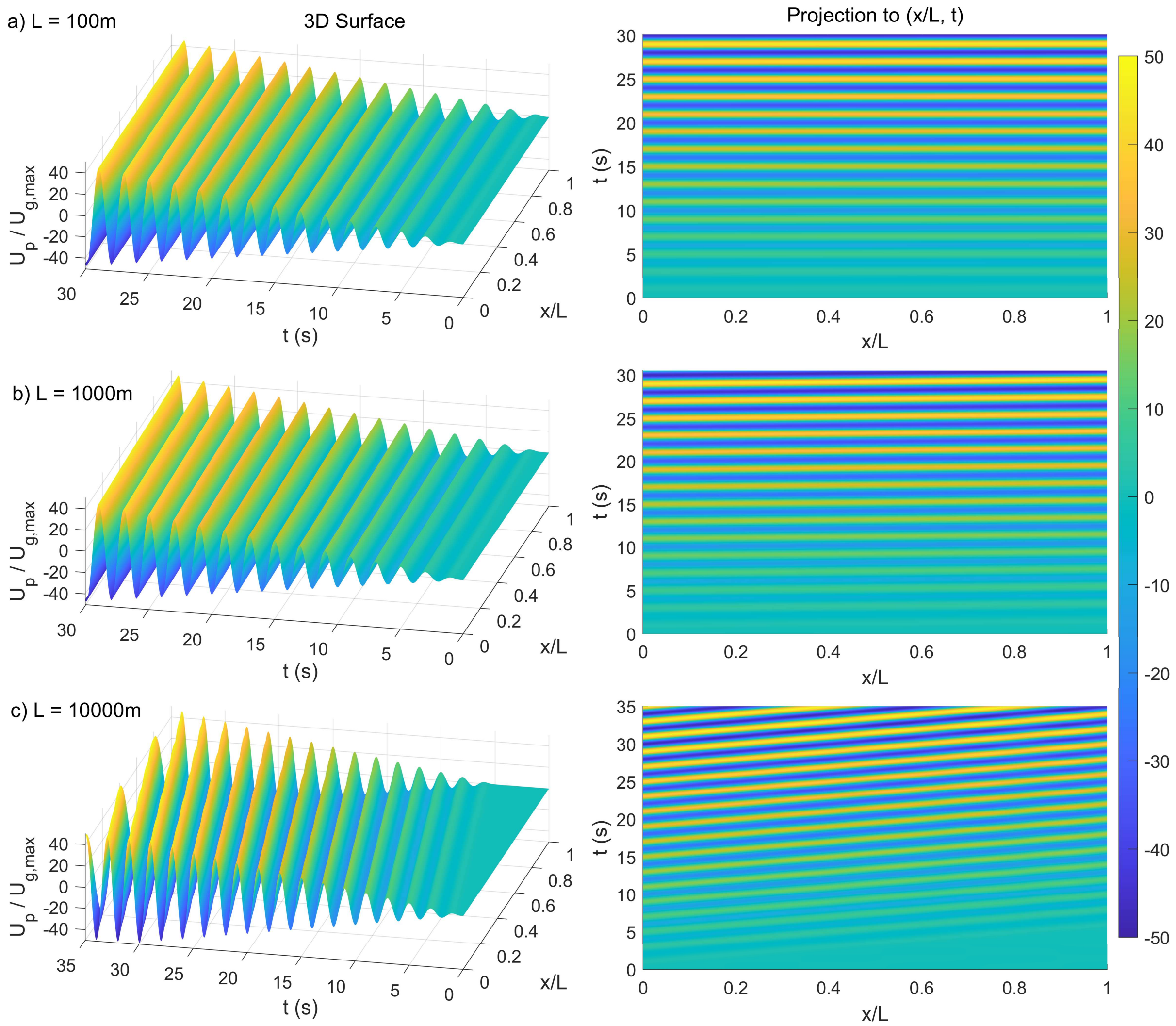

Figure 16. Variation of the relative displacement $U_p/U_{g,max}$ along the axis of the filled water pipeline ($x/L$) during the seismic wave propagation with $f$ = 0.5, considering reduced soil stiffness ($k_l$ = 11 kN/m$^2$), for various pipeline lengths: a) $L$ = 100 m; b) $L$ = 1000 m; c) $L$ = 10000 m.

Figure 17 illustrates the variation of the lateral displacement ratio $U_{p,max}$ / $U_{g,max}$ along the pipeline axis ($x/L$) at time intervals of 10s, considering a forcing frequency equal to the fundamental frequency of the unfilled ($f_1$ = 17.48 Hz) and filled water pipeline (7.68 Hz). The displacement ratio $U_{p,max}$ / $U_{g,max}$ increases with time, exceeding values of 2000 and 500 for the unfilled and filled water pipeline, respectively, for all system lengths. The soil-structure amplification at resonance increases for greater fundamental frequencies because of the greater number of cycles within the considered time history.

The predominant deformation response of the system at resonance is determined by the natural vibration mode whose average shape function wavelength is given by $\lambda = C_{ph}/f_1$, resulting in 114 m and 260 m, for the unfilled and filled water pipeline, respectively. Consequently, higher vibration modes are excited in the unfilled water pipelines, with mode numbers $n$ equal to 3, 19, and 179 corresponding to system lengths of 100 m, 1000 m, and 10000 m, respectively. These modes fall within the third part of the frequency spectrum, where the shape functions assume their maximum values at the beam ends, leading to a greater dynamic amplification in this location ($x/L$ = 1) under TGD, as shown in Figure 17.

Conversely, lower mode numbers are activated for the filled water pipelines at resonance, equal to 1, 9, and 74 for the 100m, 1000 m, and 10000 m long system respectively, resulting in more uniform amplification response along the pipeline axis (Figures 17 b, d, and f).
These elevated resonant frequencies may not fall within the range of dominant frequencies characterizing most common tectonic earthquakes, as opposed to the previous low frequency case ($f$ = 0.5 Hz). However, they are relevant to other engineering applications, including traffic induced vibration, providing further insight into the system dynamic behaviour.

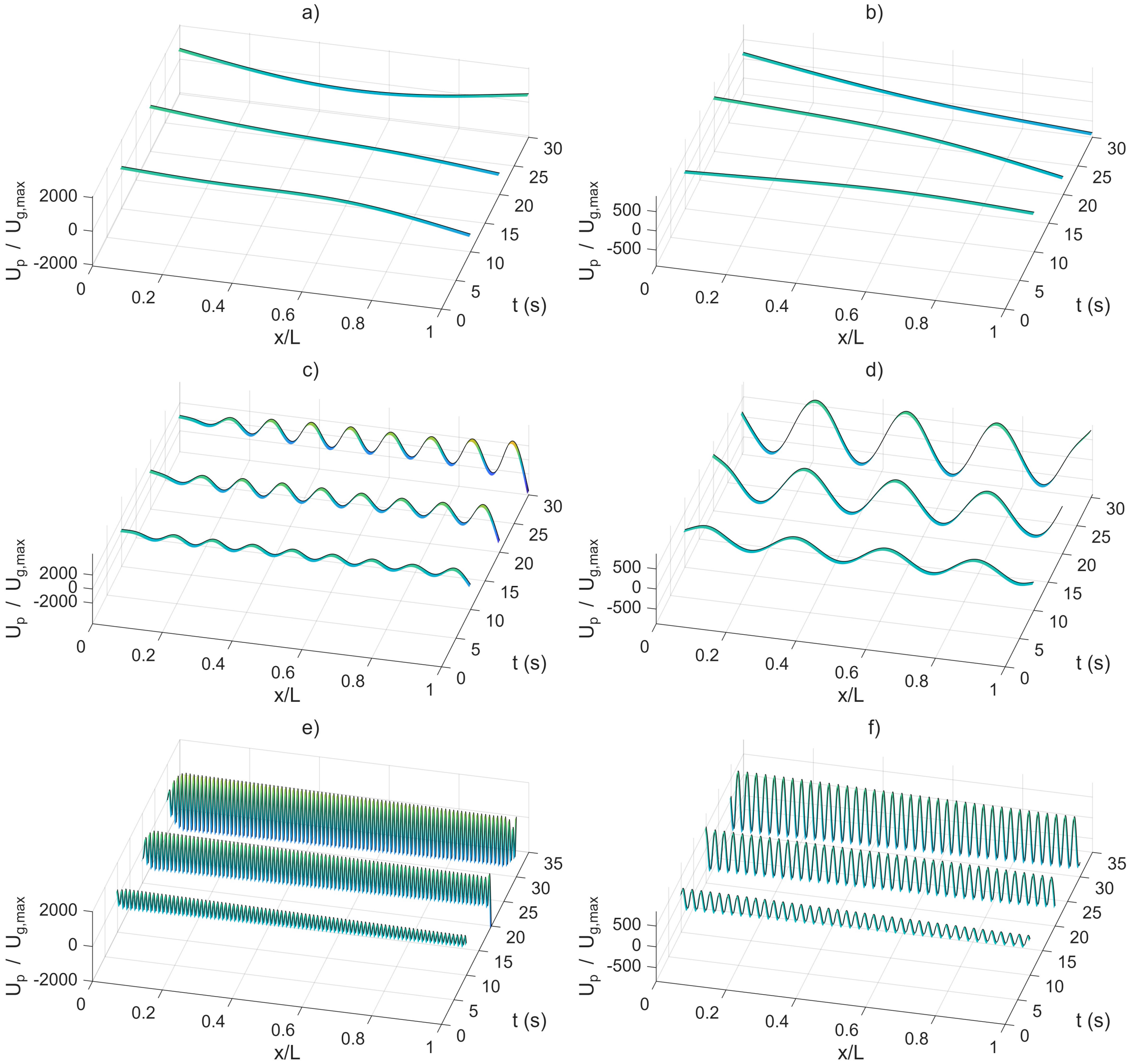

Figure 17. Variation of the relative pipeline displacement $U_p/U_{g,max}$ along its axis ($x/L$) during the seismic wave propagation with $f$ = 17.5 Hz, for various pipeline lengths and operating conditions: a) unfilled and b) filled 100 m long water pipeline; c) unfilled and d) filled 1000 m long water pipeline; e) unfilled and f) filled 10000 m long water pipeline.

## 5. Frequency response analysis

To evaluate the system response under different frequencies of the input ground motion, a frequency response analysis is performed using the developed semi-analytical model, considering different operating conditions and pipeline lengths. For comparison purposes, the duration of the harmonic ground motion is considered equal to 15 displacement cycles across the entire frequency range (0-20Hz), since the dynamic amplification depends on the number of applied loading cycles, as observed in section 4.4.

Figure 18 shows the variation of the displacement ratio ($U_{p,max}/U_{g,max}$) as a function of the ground motion frequency $f$, highlighting the resonance frequencies for the unfilled and filled water pipeline in compacted and poorly compacted backfill, for varying lengths $L$. The dynamic amplification response results less sensitive to the system length while being significantly influenced by the pipeline inertia and the soil stiffness.

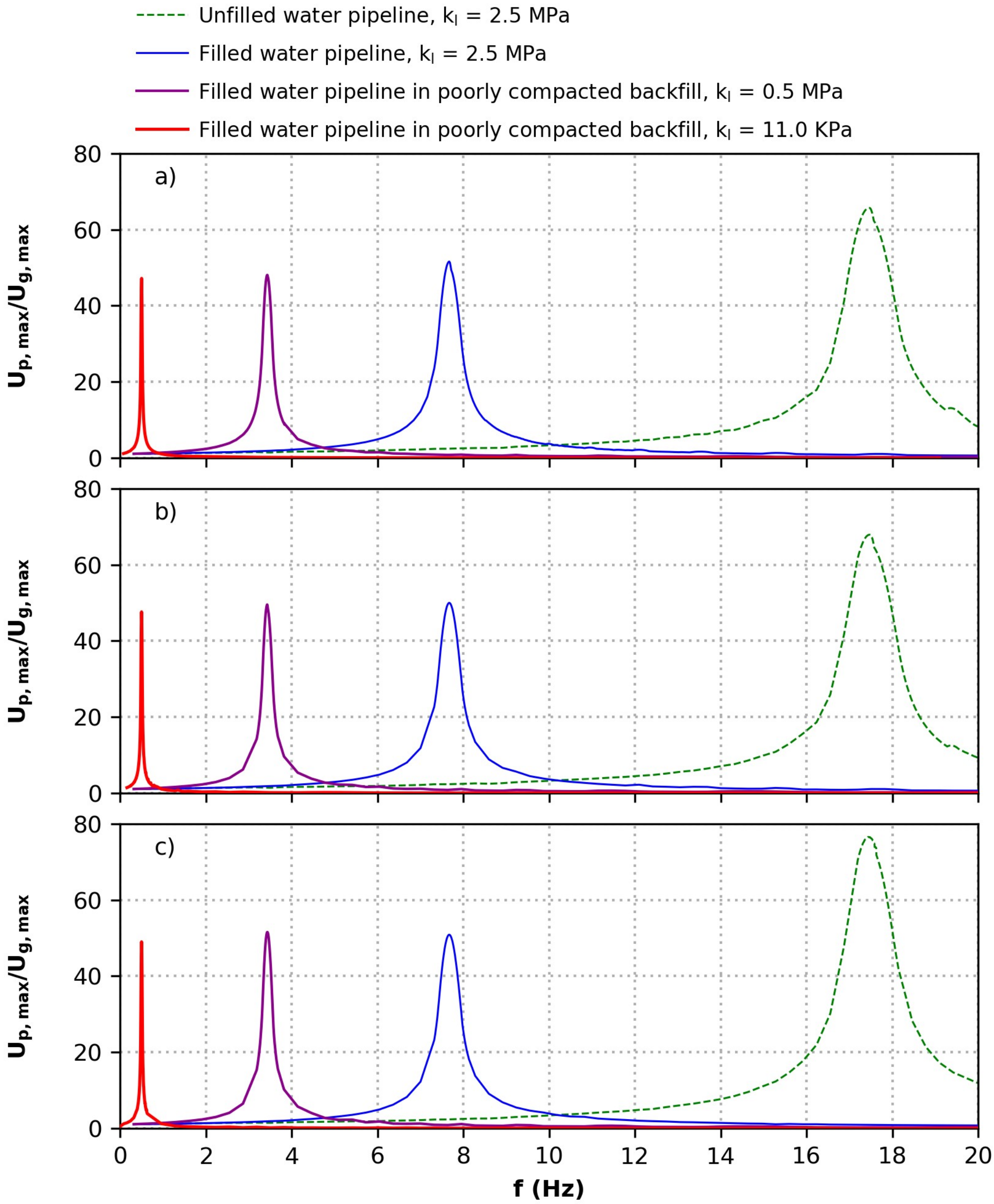

Figure 18. Maximum vibration amplitude as a function of the vibration frequency for the case of pipeline unfilled and filled with water in compacted and poorly compacted backfill, considering various pipeline lengths: a) $L$ = 100 m; b) $L$ = 1000 m; c) $L$= 10000 m.

The relative displacement $U_{p,max}/U_{g,max}$ is close to unity for low frequencies, increasing monotonically as the forcing frequency $f$ approaches the fundamental frequency of the system $f_1$, that is about 17.48 Hz and 7.68 Hz for the unfilled and filled water pipeline, respectively (Figure 5). The resonance frequency drops to 0.50 Hz and 3.43 Hz for the filled water pipeline buried in poorly compacted backfill conditions, characterized by lower values of the lateral soil stiffness, equal to $k_l$ = 11 kPa and $k_l$ = 500 kPa, respectively. These resonance frequency values may fall within the range of dominant frequencies of earthquake vibrations, making the system more vulnerable to dynamic amplification effects.

Table 6 summarizes the values of the maximum dynamic amplification $U_{p,max}/U_{g,max}$ and the corresponding frequency values, as well as the cutt-off frequencies $\tilde{\omega}_i$, for each soil-pipeline configuration. In addition, it indicates the number of modes for the first two parts of the frequency spectrum, $N_1$ and $N_2$, representing the modal density of the system in the lower frequency range ($\omega_n \leq \tilde{\omega}_2$), as well as the number $n$ of the modal shape excited at resonance. The lower modal density associated with the smaller soil stiffness results in fewer modes contributing to dynamic amplification, leading to a narrower and lower resonance bandwidth, as illustrated in Figure 18. This is more evident for the filled water pipelines, resonating at a lower modal shape number $n$ (Table 6), which results in more uniform oscillation response along the system's axis, compared to the unfilled water pipeline, as observed in Section 4.4.
Figure 19 presents the maximum dynamic amplification response $U_{p,max}/U_{g,max}$ for the unfilled and filled water pipelines ($L$ = 100 m) as a function of the soil stiffness $k_l$, and the natural frequency $\tilde{f}_2 = \tilde{\omega}_2/2\pi$. The displacement ratio $U_{p,max}/U_{g,max}$ increases monotonically with the soil stiffness, from a minimum value of about 47 for poorly compacted backfill ($k_l$ = 11 kPa) to 66 and 52 for the unfilled and filled water pipeline ($L$ = 100 m) buried in compacted soil ($k_l$ = 2.5 MPa), respectively (Figure 19 a). Consequently, the displacement ratio $U_{p,max}/U_{g,max}$ increases with the cutt-off frequency $\tilde{\omega}_2 = \sqrt{k_l/m_l}$, particularly for the unfilled water pipeline, characterized by lower inertia, as shown in Figure 19 b.
Although the resonant frequencies of the pipeline buried in compacted soil may fall outside the range of dominant frequencies of low frequency earthquakes, they are relevant for other vibration sources, including traffic, giving a better understanding of the system dynamic behaviour.

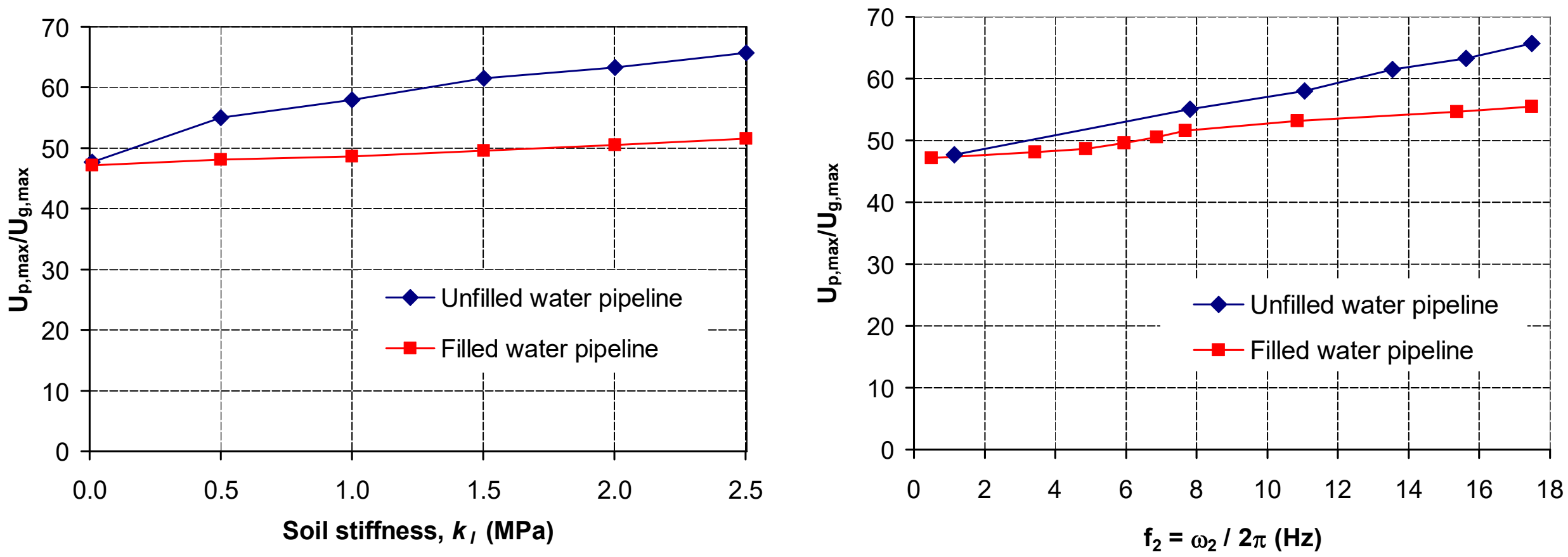

a) b)

Figure 19. Variation of the maximum dynamic amplification $U_{p,max}/U_{g,max}$ for the unfilled and filled water pipeline ($L$ = 100 m) as a function of the: a) soil stiffness $k_l$, and natural frequency $\tilde{\omega}_2 = \sqrt{k_l/m_l}$.

Table 6. Maximum dynamic amplification $U_{p,max}/U_{g,max}$ and corresponding frequency values resulting from the frequency response analysis, considering various pipeline lengths and operating conditions, characterized by the cut-off frequencies $\tilde{\omega}_i$, the number $n$ of exited modal shape at resonance, and the total number of modes for the first two parts of the frequency spectrum, $N = N_1 + N_2$, representing the modal density of the system in the low frequency range ($\omega_n \leq \tilde{\omega}_2$)

| | $L$ (m) | $U_{p,max}/U_{g,max}$ | $\omega(U_{p,max})$ (rad/s) | $f(U_{p,max})$ (Hz) | $m_l$ (kg) | $k_l$ (KPa) | $\tilde{\omega}_1$ (rad/s) | $\tilde{f}_1$ (Hz) | $\tilde{\omega}_2$ (rad/s) | $\tilde{f}_2$ (Hz) | $\tilde{\omega}_3$ (rad/s) | $\tilde{f}_3$ (Hz) | $N_1$ | $N_2$ | $N = N_1+N_2$ | $n$ |
|---|---|---|---|---|---|---|---|---|---|---|---|---|---|---|---|---|
| Unfilled water pipeline | 100 | 65.68 | 109.68 | 17.46 | | | | | | | | | 1 | 1 | 2 | 3 |
| | 1000 | 67.87 | 109.76 | 17.47 | 207.56 | 2503 | 109.8121 | 17.48 | 109.8130 | 17.48 | 5890.236 | 937.46 | 2 | 6 | 8 | 19 |
| | 10000 | 76.52 | 109.58 | 17.44 | | | | | | | | | 2 | 67 | 69 | 179 |
| Filled water pipeline | 100 | 51.60 | 48.20 | 7.67 | | | | | | | | | 1 | 1 | 2 | 1 |
| | 1000 | 49.98 | 48.20 | 7.67 | 1074.99 | 2503 | 48.2525 | 7.68 | 48.2529 | 7.68 | 2588.227 | 411.93 | 2 | 6 | 8 | 9 |
| | 10000 | 50.86 | 48.20 | 7.67 | | | | | | | | | 2 | 67 | 69 | 74 |
| Filled water pipeline | 100 | 48.055 | 21.560 | 3.43 | | | | | | | | | 1 | 1 | 1 | 1 |
| | 1000 | 49.530 | 21.540 | 3.43 | 1074.99 | 501 | 21.5793 | 3.43 | 21.5794 | 3.43 | 2588.227 | 411.93 | 2 | 3 | 5 | 4 |
| | 10000 | 51.558 | 21.560 | 3.43 | | | | | | | | | 2 | 30 | 32 | 36 |
| Filled water pipeline | 100 | 47.22 | 3.14 | 0.50 | | | | | | | | | 1 | 1 | 2 | 2 |
| | 1000 | 47.62 | 3.14 | 0.50 | 1074.99 | 11 | 3.1457 | 0.50 | 3.1457 | 0.50 | 2588.227 | 411.93 | 1 | 1 | 2 | 1 |
| | 10000 | 49.00 | 3.14 | 0.50 | | | | | | | | | 2 | 4 | 6 | 7 |

## 6. Discussion of modelling limitations

Several limitations of the present study are important to acknowledge. In current engineering design practice, soil–pipeline interaction is typically represented using uniaxial nonlinear springs in the three orthogonal directions, calculated according to widely adopted guidelines (e.g., ASCE 1984; ALA 2001, 2005). These standards, however, neglect kinematic soil–structure interaction effects and do not incorporate damping, which remains a highly uncertain parameter for underground systems. Internal structural damping within buried pipelines is generally considered small relative to soil damping, and typically neglected (Hindy and Novak, 1980; Zerva et al., 1988; Das et al., 2023). Consequently, many advanced research studies have employed undamped soil–structure interaction formulations in dynamic analyses of buried pipelines, both through analytical models (Koike, 1985; Sakurai and Takahashi, 1969; Kouretzis et al., 2007; Hosseini and Roudsari, 2015; Manolis et al., 2020) and finite-element simulations (Lee et al., 2009; Zhang et al., 2025), which further supports the modelling assumptions adopted in the present study. Nevertheless, the influence of soil and structural damping on system performance may require more detailed evaluation.

Another limitation arises from the assumption of constant soil stiffness in the Winkler foundation representation. While appropriate for moderate ground deformations and shorter system lengths, this assumption does not capture nonlinear soil behaviour, nor spatial variability in geotechnical properties, which may influence the response under large ground deformation or in heterogeneous soil conditions.

The present study also neglects the effect of internal water velocity on pipeline vibration. For fully filled buried pipelines operating under low-velocity flow, the internal fluid behaves in the limit case of a quasi-stationary medium and can be modelled as added mass, which is consistent with standard engineering practice. Although pipeline operation is inherently transient, steady-state and quasi-static models are widely used for pipeline design, flow estimation, and system optimisation (López-Benito et al., 2016), further supporting the quasi-stationary assumption adopted here. Nonetheless, for higher internal velocities or transient hydraulic phenomena, fluid–structure interaction effects may require more detailed modelling.

Overall, the proposed model is intended as a baseline analytical framework for understanding the fundamental dynamic mechanisms governing buried Timoshenko beams on elastic foundation. It provides benchmark results that can be used to calibrate more complex numerical simulations or experimental investigations. Future model extensions may incorporate system nonlinearities, viscoelastic and radiation damping, spatially variable soils, and the influence of internal fluid velocity and operating conditions.

## Conclusions

A new soil-structure interaction model is proposed for evaluating the dynamic response of Timoshenko beams on Winkler foundation subjected to transverse transient ground displacements perpendicularly to their axis. The obtained closed-form analytical solution of the governing differential equation showed that the vibration spectrum consists of four parts, separated by three transition frequencies, that can be part of the spectrum, depending on the applied boundary conditions. Across each transition frequency, the oscillatory characteristics of the vibration modes change as a function of the system's inertia and stiffness, significantly affecting the dynamic amplification response of the buried Timoshenko beam under TGD.
The accuracy of the proposed semi-analytical formulation was demonstrated through a case study of a buried steel water pipeline subjected to transverse TGD. The predicted natural frequencies, mode shapes, and forced-vibration responses showed excellent agreement with both modal and dynamic finite element analyses, confirming the correctness and reliability of the analytical solution. Additional validation was achieved through comparison with high-resolution experimental data from free-vibration tests and multi-point shaking-table experiments, further confirming the accuracy and robustness of the proposed model.

The frequency response analysis performed using the developed model showed that the ratio of the maximum pipe displacement to the maximum ground displacement, $U_{p,max} / U_{g,max}$, is close to unity for low frequencies, increasing monotonically until resonance. The latter is reached as the forcing frequency approaches the system's fundamental frequency, which depends on its mechanical and geometric properties. At resonance, the relative displacement amplitudes of the unfilled water pipeline are greater compared to the water filled pipeline, because of the greater inertia of the latter and the greater forcing frequency leading to significant dynamic amplification.

The water filled pipelines buried in poorly compacted backfill exhibited a narrower and lower resonance bandwidth, because the lower modal density associated with the smaller soil stiffness, resulting in fewer modes contributing to the dynamic amplification. In this case study, the resonance frequencies for  the water pipelines surrounded by poorly compacted soil resulted within the range of dominant frequencies of earthquake vibrations, requiring accurate seismic analysis.

Although the resonant frequencies of the pipeline buried in well compacted soil fall outside the range of dominant frequencies of low frequency earthquakes, they can be relevant for other vibration sources like high frequency seismic, traffic and railway loadings, which is a subject of the future study.

In conclusion, the proposed methodology provides a robust analytical framework for evaluating the primary factors impacting the dynamic behavior of buried beams, giving a deeper understanding of the system response under various sources of ground vibration.

## Acknowledgments

The author gratefully acknowledges Prof. Kenichi Soga for the valuable discussions and constructive feedback, and Prof. Thomas O'Rourke for his insightful guidance. The experimental data from the shaking-table tests kindly provided by Prof. Junyan Han and Prof. M. Hesham El Naggar are sincerely appreciated. The author also thanks Dr. Luis Gutiérrez for sharing the experimental data on the flexural vibration of Timoshenko beams used to validate the analytical model developed in this study.

## List of symbols

$\mathbf{A}$ = coefficient matrix for the homogenous system
$A$, $B$, $C$ = coefficients of the biquadratic equation associated with the differential equation governing the spatial function $\phi(x)$
$a_\Delta$, $b_\Delta$, $c_\Delta$ = coefficients of the biquadratic equation associated with the discriminant $\Delta$
$a_{ij}$, $b_{ij}$, $c_{ij}$, $d_{ij}$, $e_{ij}$ = elements of the coefficient matrix $\mathbf{A}$ for each part of the frequency spectrum.
$\alpha, \beta$ = eigenvalue paprameters
$A_b$ = beam cross-sectional area
$C_1$, $C_2$, $C_3$, $C_4$ = unknown constants of the general solution of the shape function $\phi(x)$
$C_n$ = generalized damping for the $n$th mode of the system
$C_{ph}$ = apparent wave propagation velocity
$\boldsymbol{C}$ = set of the complex numbers
$D$ = pipeline outer diameter
$\Delta$ = discriminant of the biquadratic equation
$dx$ = beam element length
$E$ = beam modulus of elasticity
$\phi$ = soil friction angle

$\phi_n(x)$ = modal shape corresponding to the $n$-th vibration mode
$f$ = frequency (Hz)
$\tilde{f}_c$ = cutt-off frequency predicted by the Timoshenko beam theory, representing the transition values between two different solutions of the differential equation governing the spatial function $\phi(x)$
$\tilde{f}_i$ = cutt-off frequencies (Hz), representing the transition values between two different solutions of the differential equation governing the spatial function $\phi(x)$
$f_r$ = frictional force per unit length of pipe
$\gamma$ = soil density
$G$ = Elastic shear modulus
$H$ = burial depth to pipe springline
$J$ = beam second moment of inertia
$k$ = Timoshenko beam shear coefficient
$k_l$ = the stiffness of the lateral soil spring
$K_0$ = Coefficient of lateral soil pressure at rest
$K_M$ = Coefficient for the moment boundary condition at the beam free ends
$K_n$ = generalized stiffness for the $n$th mode of the system
$K_T$ = Coefficient for the shear boundary condition at the beam free ends
$L$ = length of the buried beam/pipeline
$l_c = \sqrt[4]{4EJ/k_l}$ characteristic length of the system
$\lambda$ = variable of the characteristic equation associated with the differential equation governing the spatial function $\phi(x)$
$\lambda_i$ = solutions of the characteristic equation associated with the differential equation governing the spatial function $\phi(x)$
$M$ = beam bending moment
$M_n$ = generalized mass for the $n$th mode of the system
$\mu$ = coefficient of friction at the soil-structure interface
$m_l$ = linear mass of the beam per unit length
$n$ = mode vibration number
$N_i$ = number of vibration modes within the $i$-th part of the frequency spectrum
$\nu$ = beam's Poisson's ratio
$p_u$ = lateral soil reaction per unit pipe length
$P_n(t)$ = generalized force for the $n$th mode of the system
$Q$ = beam shear force
$\varphi$ = beam cross-sectional rotation
$q_n(t)$ = modal coordinate of the $n$-th vibration mode
$r$ = radius of gyration of the beam cross section
$t$ = pipeline thickness
$U_g$ = field of ground displacement
$U_{g,max}$ = maximum ground displacement
$u_l$ = lateral soil-spring maximum elastic deformation
$u_0$ = relative soil-pipe displacement at the onset of friction sliding
$U_p$ = pipeline displacement
$U_{p,max}$ = maximum pipeline displacement
$\tilde{\omega}_i$ = cutt-off angular frequencies (rad/s), representing the transition values between two different solutions of the differential equation governing the spatial function $\phi(x)$
$\omega_n$ = vibration frequency of the $n$-th mode
$x$ = horizontal distance along the pipeline system
$\xi_n$ = damping ratio for the $n$th mode of the system
X = unknown vector for the homogenous system
$y_G$ = transient ground displacement
$y$ = transverse displacement of the beam centroid
$y_n(x,t)$ = displacement corresponding to the $n$-th vibration modes